\documentclass[11pt]{article}

\usepackage{amssymb,amsmath,latexsym}

\begin{document}

\newtheorem{thm}{Theorem}[section]
\newtheorem{lemma}[thm]{Lemma}
\newtheorem{defn}[thm]{Definition}
\newtheorem{prop}[thm]{Proposition}
\newtheorem{corollary}[thm]{Corollary}
\newtheorem{remark}[thm]{Remark}
\newtheorem{example}[thm]{Example}

\numberwithin{equation}{section}

\def\ee{\varepsilon}
\def\qed{{\hfill $\Box$ \bigskip}}
\def\MM{{\cal M}}
\def\BB{{\cal B}}
\def\LL{{\cal L}}
\def\FF{{\cal F}}
\def\EE{{\cal E}}
\def\QQ{{\cal Q}}
\def\AA{{\cal A}}
\def\CC{{\cal C}}
\def\R{{\bf R}}
\def\N{{\mathbb N}}
\def\E{{\bf E}}
\def\F{{\bf F}}
\def\H{{\bf H}}
\def\P{{\bf P}}
\def\Q{{\bf Q}}
\def\S{{\bf S}}
\def\J{{\bf J}}
\def\K{{\bf K}}
\def\F{{\bf F}}
\def\A{{\bf A}}
\def\loc{{\bf loc}}
\def\eps{\varepsilon}
\def\semi{{\bf semi}}
\def\wh{\widehat}
\def\pf{\noindent{\bf Proof.} }
\def\dim{{\rm dim}}

\title{\Large \bf On dual processes
of non-symmetric diffusions with measure-valued drifts}
\author{Panki Kim\\
Department of Mathematics\\
Seoul National University\\
Seoul 151-742, Republic of Korea\\
Email: pkim@snu.ac.kr \smallskip \\
URL: www.math.snu.ac.kr/$\sim$pkim\\
Telephone number: 82-2-880-4077\\
Fax number: 82-2-887-4694 \\
 and
\smallskip \\
Renming Song\thanks{The research of this author is supported in part
by a joint
US-Croatia grant INT 0302167.}\\
Department of Mathematics\\
University of Illinois \\
Urbana, IL 61801, USA\\
Email: rsong@math.uiuc.edu\\
URL:www.math.uiuc.edu/$\sim$rsong\\
 Telephone number: 1-217-244-6604\\
Fax number: 1-217-333-9576  }
\maketitle

\bigskip

\begin{abstract}
In this paper, we study properties of the dual process and
Schr\"{o}dinger-type operators
of a  non-symmetric diffusion with measure-valued drift.
Let $\mu=(\mu^1, \dots, \mu^d)$ be such that each $\mu^i$ is
a signed measure
on $\R^d$ belonging to the Kato class $\K_{d, 1}$.
A diffusion with drift $\mu$ is a
diffusion process in $\R^d$ whose generator can be informally
written as $L+\mu\cdot\nabla$ where
$L$ is a uniformly elliptic
differential operator.
When each $\mu^i$ is given by $U^i(x)dx$ for some function
$U^i$, a diffusion with drift $\mu$ is a diffusion in $\R^d$
with generator $L+U\cdot\nabla$.
In \cite{KS, KS1}, we have already studied properties of
diffusions with
measure-value drifts in bounded domains.
In this paper we discuss the potential theory of
the  dual process and Schr\"{o}dinger-type operators
of a diffusion with measure-valued drift.

We show that a killed diffusion process with measure-valued drift
in any bounded domain has a dual process with respect to a
certain reference measure.
For an arbitrary bounded domain,
we show that a scale invariant
Harnack inequality is true for the dual process.
We also show that, if the domain is  bounded $C^{1,1}$, the
boundary Harnack principle for the dual process is true and
the (minimal) Martin boundary for the dual process can be identified
with the Euclidean boundary. It is also shown that the harmonic measure for
the dual process is locally comparable to that of the
$h$-conditioned Brownian motion with
$h$ being the ground state.

Informally the Schr\"{o}dinger operator we consider is
 $L+\mu\cdot\nabla + \nu$ where $L$ is uniformly elliptic,
$\mu$ is a vector-valued signed measure
 in $\R^d$ and $\nu$ is a signed measure in $\R^d$.
Under the gaugeability assumption,
if the domain is
 bounded Lipschitz, the (minimal) Martin boundary
for  the Schr\"odinger operator obtained from the diffusion
with measure-value drift can
be identified with the Euclidean boundary.
\end{abstract}

\vspace{.6truein}

\bigskip


\bigskip

\begin{doublespace}

\section{Introduction}

In this paper, we continue discussing diffusions with measure-valued drift
from \cite{ KS, KS1}.

Throughout this paper, we always assume that $d\ge 3$. First we recall
the definition of the Kato class $\K_{d, \alpha}$
for $\alpha\in (0, 2]$.
For any function $f$ on $\R^d$ and $r>0$, we define
$$
M^ \alpha_f(r)=\sup_{x\in \R^d}\int_{|x-y|\le r}\frac{|f|(y)dy}
{|x-y|^{d- \alpha}}, \quad \quad 0 <  \alpha \le 2.
$$
In this paper, we mean, by a signed measure, the difference of two
nonnegative measures at most one of which can have infinite
total mass.
For any signed measure $\nu$ on $\R^d$, we use
$\nu^+$ and $\nu^-$ to denote its positive and negative parts, and
$|\nu|=\nu^++\nu^-$ its total variation.
For any signed measure $\nu$  on $\R^d$ and any $r>0$, we define
$$
M^\alpha_{\nu}(r)=\sup_{x\in \R^d}\int_{|x-y|\le r}\frac{|\nu|(dy)}
{|x-y|^{d- \alpha}}, \quad  \quad 0 < \alpha  \le 2.
$$

\medskip

\begin{defn}\label{d:kc} Let $0 <  \alpha \le 2.$
We say that a function $f$ on $\R^d$ belongs to the Kato
class $\K_{d, \alpha}$ if
$\lim_{r\downarrow 0}M^\alpha_{f}(r)=0$.
We say that a signed Radon measure $\nu$ on $\R^d$
belongs to the Kato class $\K_{d, \alpha}$ if
$\lim_{r\downarrow 0}M^\alpha_{\nu}(r)=0$.
We say that a $d$-dimensional vector valued function $V =
(V^1, \cdots , V^d)$ on $\R^d$ belongs to the Kato
class $\K_{d, \alpha}$ if  each $V^i$  belongs to the Kato
class $\K_{d, \alpha}$.
We say that a $d$-dimensional vector valued signed Radon measure
$\mu = (\mu^1, \cdots, \mu^d)$ on $\R^d$
belongs to the Kato class $\K_{d, \alpha}$ if each $\mu^i$
belongs to the Kato
class $\K_{d, \alpha}$.
\end{defn}

\medskip

Rigorously speaking a function $f$ in $\K_{d, \alpha}$
may not give rise to a signed measure $\nu$ in $\K_{d, \alpha}$
since it may not give rise to a signed measure at all. However, for
the sake of simplicity we use the convention that whenever we write
that a signed measure $\nu$ belongs to $\K_{d, \alpha}$
we are implicitly assuming that we are covering the case of
all the functions
in $\K_{d, \alpha}$ as well.

Throughout this paper we assume that $\mu=(\mu^1, \dots, \mu^d)$ is
fixed with each $\mu^i$ being a signed measure on $\R^d$ belonging
to $\K_{d, 1}$.  We also assume that the operator
$L$ is either $L_1$ or $L_2$ where
$$
L_1:= \frac12\sum_{i,j=1}^{d} \partial_i (a_{ij} \partial_j)
\qquad \mbox{and} \qquad
L_2:=\frac12 \sum_{i,j=1}^{d} a_{ij} \partial_i \partial_j
$$
with $\A:=(a_{ij})$ being $C^{1}$ and uniformly elliptic but not necessarily
symmetric.

Informally, when $a_{ij}$ is symmetric,
a diffusion process in $\R^d$ with drift $\mu$ is a
diffusion process
in $\R^d$ with generator
$ L+\mu\cdot\nabla$.
When each $\mu^i$ is given by $U^i(x)dx$ for some function
$U^i$, a diffusion process with drift $\mu$ is a diffusion in $\R^d$
with generator $L +U\cdot\nabla$
and it is a solution to the SDE
$
dX_t=dY_t+U(X_t)\cdot dt
$
where $Y$ is  a diffusion in $\R^d$
with generator $L$

To give the precise definition of a diffusion with drift $\mu$
in $\K_{d, 1}$,
we fix a non-negative smooth radial function $\varphi(x)$ in
$\R^d$ with supp$[\varphi] \subset B(0,1)$
and $\int \varphi(x) dx =1 $. For any positive integer $n$, we
put $\varphi_n(x)=2^{nd}\varphi(2^nx)$.
For $1\le i\le d$, define
$$
U^i_n(x)= \int \varphi_n(x-y) \mu^i(dy).
$$
Put $U_n(x)=(U^1_n(x), \dots, U^d_n(x))$.

In the definition below, we assume $a_{ij}$ is symmetric.

\medskip

\begin{defn}\label{d:sde2}
Suppose
$\mu=(\mu^1, \dots, \mu^d)$ is such that each $\mu^i$ is a signed measure on
$\R^d$ belonging to the Kato class $\K_{d, 1}$.
A diffusion with drift $\mu$ is a family of probability measures
$\{\P_x: x\in \R^d\}$ on $C([0, \infty), \R^d)$, the space
of continuous $\R^d$-valued functions on $[0, \infty)$, such that
under each $\P_x$ we have
$$
X_t=x+Y_t +A_t
$$
where
\begin{description}
\item{(a)} $A_t=\lim_{n\to \infty}\int^t_0U_n(X_s)ds$ uniformly in $t$
over finite intervals, where the convergence is in probability;
\item{(b)} there exists a subsequence $\{n_k\}$ such that
$$
\sup_k\int^t_0|U_{n_k}(X_s)|ds<\infty
$$
almost surely for each $t>0$;
\item{(c)} $Y_t$ is a diffusion in $\R^d$
starting from the origin with generator $L$.
\end{description}
\end{defn}

\medskip

The existence and uniqueness
of $X$ were established in \cite{BC} (see Remark 6.1 in \cite{BC}).
For any open set $D$, we use $\tau_D$ to denote the first exit time of $D$,
i.e., $\tau_D=\inf\{t>0: \, X_t\notin D\}$.
Given  an open set $D\subset \R^d$, we define
$X^D_t(\omega)=X_t(\omega)$ if $t< \tau_D(\omega)$
and $X^D_t(\omega)=\partial$ if $t\geq  \tau_D(\omega)$,
where $\partial$ is a cemetery state. The process
$X^D$ is called a killed diffusion with drift $\mu$ in $D$.
In this paper we discuss properties of $X^D$ when $D$ is a
bounded domain.

When $a_{ij}$ is not symmetric, we use a simple
reduction;
Let $Y_t$ be the diffusion in $\R^d$
with generator
$$
\frac14\sum_{i,j=1}^{d} \partial_i ((a_{ij}+a_{ji}) \partial_j).
$$
Note that
$$
\sum_{i,j=1}^{d} a_{ij} \partial_i \partial_j
= \sum_{i,j=1}^{d}\frac12(a_{ij}+a_{ji}) \partial_i \partial_j=
\sum_{i,j=1}^{d} \frac12\partial_i ((a_{ij}+a_{ji}) \partial_j)-
 \sum_{i,j=1}^{d}  \frac12\partial_i (a_{ij}+a_{ji}) \partial_j
$$
and
$$\sum_{i,j=1}^{d} \partial_i (a_{ij} \partial_j)
= \sum_{i,j=1}^{d} a_{ij} \partial_i \partial_j +
 \sum_{i,j=1}^{d} (\partial_i a_{ij}) \partial_j
=\sum_{i,j=1}^{d}\frac12 \partial_i ((a_{ij}+a_{ji}) \partial_j)
+\sum_{i,j=1}^{d}\frac12  \partial_i (a_{ij}-a_{ji}) \partial_j.
$$
Since, for  any bounded domain $D$,
$$
\left( \sum_{i=1}^{d}\frac14  \partial_i  (a_{i1}+a_{1i})|_D,
\cdots,  \sum_{i=1}^{d}\frac14  \partial_i  (a_{id}+a_{di})|_D\right)
 \quad  \mbox{and}\quad
\left(\sum_{i=1}^{d}\frac14  \partial_i (a_{i1}-a_{1i})|_D, \cdots,
\sum_{i=1}^{d}\frac14  \partial_i (a_{id}-a_{1d})|_D \right)
$$
are in $\K_{d,1}$,
we construct $X_t$ with a drift which is either
$$
\left( \mu^1+\sum_{i=1}^{d}\frac14  \partial_i (a_{i1}+a_{1i})|_Ddx,
\cdots, \mu^d +\sum_{i=1}^{d}\frac14  \partial_i (a_{id}+a_{di})|_Ddx\right)
$$
or
$$
\left(\mu^1+\sum_{i=1}^{d}\frac14  \partial_i (a_{i1}-a_{1i})|_Ddx, \cdots,
 \mu^d+\sum_{i=1}^{d}\frac14  \partial_i (a_{id}-a_{1d})|_Ddx \right)
$$ as in the Definition \ref{d:sde2}.
Then the generator of the killed diffusion process
$X^D$ in $D$ can be informally written as
$L+ \mu \cdot \nabla$ where $L$ is either $$
\frac12\sum_{i,j=1}^{d} \partial_i (a_{ij} \partial_j)
\qquad \mbox{and} \qquad
\frac12 \sum_{i,j=1}^{d} a_{ij} \partial_i \partial_j.
$$
with $\A:=(a_{ij})$ being $C^{1}$ and uniformly elliptic but not necessarily
symmetric.

Throughout this paper we assume that $X^D$ is the process constructed
above.
In \cite{KS} (also see section 6 in \cite{KS1}), we showed that $X$
has a density $q(t, x, y)$ which is continuous on
$(0, \infty)\times \R^d\times \R^d$ and that there exist positive
constants $c_i$, $i=1, \cdots, 9$, such that
\begin{equation}\label{e:2sge}
c_1e^{-c_2 t}t^{-\frac{d}{2}}
e^{-\frac{c_3 |x-y|^2}{2t}}
\le q(t, x, y)\le c_4
e^{c_5 t}t^{-\frac{d}{2}}e^{-\frac{c_6|x-y|^2}{2t}}
\end{equation}
and
\begin{equation}\label{e:gedbmsd}
|\nabla_x q(t, x, y)|\le C_7
e^{c_8 t}t^{-\frac{d+1}{2}}e^{-\frac{c_9|x-y|^2}{2t}}
\end{equation}
for all $(t, x, y)\in (0, \infty)\times \R^d\times \R^d$.
We also showed that,
for every bounded $C^{1,1}$ domain $D$ (see below for the definition),
$X^D$ has a  density $q^{D}$ which is continuous
on $(0, \infty)\times D\times D$ and
that for any $T>0$,
there exist positive constants $c_i$, $i=10, \dots, 14,$ such that
\begin{equation}\label{e:dekbmsd}
c_{10}
t^{-\frac{d}2}(1\wedge \frac{\rho(x)}{\sqrt t})
(1\wedge \frac{\rho(y)}{\sqrt t})
e^{-\frac{c_{11}|x-y|^2}t}
\le q^{D}(t, x, y)\le
c_{12}(1\wedge \frac{\rho(x)}{\sqrt t})
(1\wedge \frac{\rho(y)}{\sqrt t})
t^{-\frac{d}2}e^{-\frac{c_{13}|x-y|^2}{t}}
\end{equation}
and
\begin{equation}\label{e:gedkbmsd}
 |\nabla_x q^{D}(t, x, y)|\le
c_{14}(1\wedge \frac{\rho(y)}{\sqrt t})
t^{-\frac{d+1}2}e^{-\frac{c_{13}|x-y|^2}{t}}
\end{equation}
for all $(t, x, y)\in (0, T]\times D\times D$,
where  $a \wedge b := \min \{ a, b\}$, $\rho(x)$ is
the distance between $x$ and
$\partial D$.
Using the estimates above we studied the potential theory of
$X$ in \cite{KS} and \cite{KS1}. More precisely, we proved
the boundary Harnack principle holds for nonnegative harmonic functions
of $X$ in bounded Lipschitz domains and identified the Martin
boundary of the killed process $X^D$ when $D$ is a bounded
Lipschitz domain.

In general, the process $X$ does not have a dual and therefore
results for Markov processes under the duality hypotheses,
like the the general conditional gauge theorems of \cite{C2} and \cite{CS6}
or the stability of Martin boundary under perturbation of \cite{CK2},
can not be applied to $X$ directly. In this paper we will prove
that, for any bounded domain $D$,  $X^D$ has a dual process
with respect to a certain reference
measure and the dual process is a continuous Hunt process
on $D$ with the strong Feller property.
By using this duality, we can apply the general conditional gauge
theorems of \cite{C2} and \cite{CS6}
and the stability result of \cite{CK2} to the present situation.

After establishing the existence of the dual process, we study
properties of the dual process. We show that a scale invariant
Harnack inequality is true for the dual process.
We also show that, if the domain is  bounded $C^{1,1}$, the
boundary Harnack principle for the dual process is true and
the (minimal) Martin boundary for the dual process can be identified
with the Euclidean boundary. One of the interesting fact is that
 the harmonic measure for
the dual process is locally comparable to that of the
$h$-conditioned Brownian motion with
$h$ being the ground state.

In \cite{KS4} the concept of intrinsic ultracontractivity was extended
to non-symmetric semigroups and it was proved there that
the semigroup
of a killed diffusion process in a bounded Lipschitz domain is intrinsic
ultracontractive if the coefficients of the generator of the diffusion
process are smooth.
In \cite{KS5} we will use the duality proved in this paper
to show that the Schr\"{o}dinger semigroup of
the killed process $X^D$ is intrinsic ultracontractive under very
weak assumptions on $D$.

The content of this paper is organized as follows. In Section 2,
we present some preliminary properties of the killed process $X^D$
in an arbitrary bounded domain $D$; the existence of
the dual process of $X^D$
is proved in Section 3; Section 4 contains a result on the
comparison of harmonic
measures and a scale invariant Harnack inequality for the dual process
which is used in Sections 5 and 6 to study the Martin boundary
of the dual process; and in the last section we specialize the
general conditional gauge theorems of \cite{C2} and \cite{CS6}
to the present setting
and then, using the stability result of \cite{CK2},
get the stability of Martin boundaries of $X^D$ and its dual
under perturbations.

Throughout this paper, we use the notation $a \wedge b := \min \{ a, b\}$ and
$a \vee b := \max \{ a, b\}$.
We will use the convention $f(\partial)=0$ for any function $f$ on $D$.
In this paper we will also use the following convention:
the values of the constants $c_1, c_2, \cdots$
might change from one appearance to another. The labeling
of the constants $c_1, c_2, \cdots$ starts anew in the
statement of each result. In this paper, we use ``$:=$" to denote a
definition, which is  read as ``is defined to be".

\section{Diffusion with measure-valued drift in bounded domains}

In this section we assume that
$D$ is an arbitrary bounded domain and we will discuss some basic
properties of $X^D$
that we  will  need later.

It is shown in \cite{KS1} that $X^D$ has a
 jointly continuous and strictly positive transition density
function $q^D(t,x,y)$.
Using the continuity $q^D(t,x,y)$ and the estimates (\ref{e:2sge}),
the proof of the next proposition is easy. We omit the proof.

\medskip

\begin{prop}\label{p:HC}
$X^D$ is a Hunt processes and has the strong Feller property. i.e, for every
$f \in L^{\infty}(D)$, $P^D_tf(x):=\E_x[f(X^D_t)]$ is bounded
and continuous in $D$.
\end{prop}

\medskip

We know from \cite{KS1} that
there exist positive constants $c_1$ and $c_2$ depending on $D$
via its diameter such that
for any $(t, x, y)\in (0, \infty)\times D\times D$,

\begin{equation}\label{est:6.1}
\,q^{D}(t, x, y)\,\le\, c_1
t^{-\frac{d}{2}}      e^{-\frac{c_2|x-y|^2}{2t}}.
\end{equation}
Let $G_D(x, y)$ be the Green function of $X^D$, i.e.,
$$
G_D(x, y):=\int^{\infty}_0 q^D(t, x, y)dt.
$$  $G_D(x, y)$  is finite
 for $x\neq y$ and
\begin{equation}\label{G_bd}
G_D(x, y) \, \le \,
  \frac{c_3}{|x-y|^{d-2}}
\end{equation}
for some $c_{3}=c_3
($diam$(D))>0$.
Now define
$$
h_D(x)\,:=\, \int_D G_D(y,x) dy\quad\mbox{and}\quad \xi_D(dx):=h_D(x) dx.
$$

The following result says that $\xi_D$ is a reference measure for $X^D$.

\medskip

\begin{prop}\label{e_m} For any bounded domain $D$,
 $\xi_D$ is  an excessive measure with respect to $X^D$,
i.e., for every Borel function $f \ge 0$,
$$
\int_D f(x) \xi_D(dx) \ge \int_D  \E_x\left[f(X^D_{t})\right] \xi_D(dx).
$$
Moreover, $h_D$ is a strictly positive, bounded continuous function on $D$.
\end{prop}
\pf
By the Markov property, we have for any Borel function $f \ge 0$,
\begin{eqnarray*}
 &&\int_D\E_y\left[f(X^D_{t})\right] G_D(x,y) dy
\,=\,\E_x \int_0^\infty\E_{X^D_s}\left[f(X^D_{t})\right] ds\\
&&= \int_0^\infty\E_{x}\left[f(X^D_{t+s})\right] ds
\, \le \,  \int_D f(y) G_D(x,y) dy, \quad x \in D.
\end{eqnarray*}
Integrating with respect to  $x$, we get by Fubini's theorem,
 $$\int_D\E_y\left[f(X^D_{t})\right] h_D(y) dy
\le \int_D f(y)  h_D(y) dy$$
The second claim follows from the continuity of $G_D$ and (\ref{G_bd}).
\qed

\medskip

We define a new transition density function by
$$
\overline{q}^D(t,x,y)\,:=\, \frac{ q^D(t,x,y)}{h_D(y)}.
$$
Let
$$
 \overline{G}_D(x,y)\,:=\,
\int^{\infty}_0 \overline{q}^D(t, x, y)dt\,=\,\frac{ G_D(x,y)}{h_D(y)}.
$$
Then $\overline{G}_D(x,y)$ is
the Green function of $X^D$ with respect to the reference measure $\xi_D$.

Before we discuss properties of  $ \overline{G}_D(x,y)$,
let's first recall some definitions.

\medskip

\begin{defn}\label{d:har}
Suppose $U$ is an open subset of $D$.
A Borel function $u$ defined on $U$  is said to be
\begin{description}
\item{(1)} harmonic with respect to $X^D$ in $U$ if
\begin{equation}\label{AV}
u(x)\,= \,\E_x\left[u(X^D_{\tau_{B}})\right],
\qquad x\in B,
\end{equation}
for every bounded open set $B$ with $\overline{B}\subset U$;

\item{(2)} superharmonic
with respect to $X^D$ if
$u$ is non-negative and
$$
u(x)\,\geq\,  \E_x\left[u(X^D_{\tau_{B}})\right], \qquad x\in B,
$$
for every bounded open set $B$ with $\overline{B}\subset D$;

\item{(3)} excessive for $X^D$  if $u$ is non-negative and
$$
u(x)\,\geq\,  \E_x\left[u(X^D_{t})\right]\quad
 { and }\quad u(x)\,=\, \lim_{t \downarrow 0}
\E_x\left[u(X^D_t)\right], \qquad t > 0, x\in D;
$$

\item{(4)} a potential for  $X^D$ if
it is excessive for $X^D$ and for every sequence $\{U_n\}_{n \ge 1}$
of open sets with
$\overline{U_n} \subset U_{n+1}$ and  $\cup_n U_n = D$,
$$
 \lim_{n \to \infty}
 \E_x\left[u(X^D_{\tau_{U_n}})\right]\,=\, 0;
\qquad \xi_D\mbox {-a.e. } x\in D.
$$
\item{(5)}  a pure potential for $X^D$ if
it is a potential for  $X^D$ and $$
 \lim_{t \to \infty}
 \E_x\left[u(X^D_t)\right]\,=\, 0, \qquad \xi_D\mbox {-a.e. } x\in D.
$$
\end{description}
A Borel function $u$ defined on $\overline{U}$
is said to be regular harmonic with respect to $X^D$ in $U$ if
$u$ is harmonic with respect to $X^D$ in $U$ and
(\ref{AV}) is true for $B=U$;

A Borel function $u$ defined on
$D$ is said to be
harmonic with respect to $X^D$  if it is harmonic with
respect to $X^D$ in $D$.
\end{defn}

\medskip

Since $X^D$ is a Hunt processes with the strong Feller property,
it is easy to check
that $u$ is excessive for $X^D$ if and only if $f$ is
lower semi-continuous in $D$ and superharmonic with respect to $X^D$.
(See Theorem 4.5.3 in \cite{CW} for the Brownian motion case, and the
proof there can adapted easily to the present case.)

We list some properties of the Green function $\overline{G}_D(x,y)$
of $X^D$ that  we will need later.
\begin{description}
\item{(A1)} $\overline{G}_D(x,y) > 0$ for all $(x,y) \in D \times D$;
$\overline{G}_D(x,y) =\infty$
  if and only if $x=y \in D$;
\item{(A2)} For every $x \in D$, $\overline{G}_D(x,\,\cdot\,)$ and
$\overline{G}_D(\,\cdot\,,x)$ are extended
continuous in $D$;
\item{(A3)} For every compact subset $K$ of $D$,
$\int_K \overline{G}_D(x,y) \xi_D(dy) < \infty$.
\end{description}

The above properties can be checked easily from
Theorem 2.6 in \cite{KS1} and our Proposition \ref{e_m}
and (\ref{G_bd}) above.
Thus $X^D$ is a transient diffusion
satisfying  the conditions in \cite{CR} and \cite{P}.
Applying  Theorem 1 in \cite{P}, we have that

\begin{description}
\item{(A4)}
for each $y$,  $ x \to \overline{G}_D(x,y)$ is excessive for
$X^D$ and harmonic for $X^D$ in
$D \setminus \{y\}$.
Moreover, for every
open subset $U$ of $D$, we have
\begin{equation}\label{e:GH}
\E_x[\overline{G}_D(X^D_{T_{U}},y)]=\overline{G}
_D(x,y), \qquad (x,y) \in D \times U
\end{equation}
where  $T_U:=\inf \{ t > 0: X^D_t \in U\}$. In particular,
for every $y \in D$ and $\eps > 0$,
$\overline{G}_D(\,\cdot\,, y)$
is regular harmonic in $D \setminus B(y, \eps)$
with respect to $X^D$.
\end{description}

Since the set $Z$ defined in \cite{CR}
(equation (12) on page 179) is empty, from
Corollary 2 and Theorems 5-6 in \cite{CR}, we have

\medskip

\begin{thm}\label{t:Riesz}
\begin{description}
\item{(1)}  If $u$ is a potential for $X^D$, then
there exists
a unique Radon measure $\nu$ on $D$ such that
$$
u(x) \,=\,G_D \nu (x)\, := \int_D\overline{G}
_D(x,y) \nu (dy) $$
\item{(2)} If $f$ is an excessive function for $X^D$
and $f$ is not identically zero,
then there exists
a unique Radon measure $\nu$ on $D$ and a nonnegative harmonic function $h$
for $X^D$ such that
$ f=G_D \nu + h.$
\end{description}
\end{thm}

\medskip

For $y\in D$, let $X^{D,y}$ denote the $h$-conditioned process obtained
from  $X^D$ with
$h(\cdot)= \overline{G}_D(\cdot , y)$ and
let $\E_x^y$ denote the expectation
for $X^{D,y}$ starting from $x\in D$.

\medskip

\begin{thm}\label{t:p}
For each $y$,  $ x \to \overline{G}_D(x,y)$ is a pure potential for $X^D$.
In fact, for every sequence $\{U_n\}_{n \ge 1}$
of open sets with
$\overline{U_n} \subset U_{n+1}$ and  $\cup_n U_n = D$,
$$
\lim_{n \to \infty}
\E_x\left[\overline{G}_D(X^D_{\tau_{U_n}},y)\right]\,=\, 0,
\qquad x \not= y.
$$
Moreover, for every $x,y \in D$, we have
$$
\lim_{t \to \infty} \E_x\left[\,\overline{G}_D(X^D_t, y)\right] \,=\, 0.
$$
\end{thm}

\pf Let $x \not= y \in D$.
First we see that, from (A1)-(A2), the condition (H) in \cite{LZ} holds.
Also the strict positivity of $\overline{G}_D$ and (A4) imply that
the set $W$ on page 5 in \cite{LZ} (also the set $Z$ defined in \cite{CR})
is empty.
Thus by Theorem 2 in \cite{LZ}, the lifetime $\zeta^y$ for $X^{D,y}$ is
finite $\P_x^y$-a.s. and
\begin{equation}\label{e_lim}
\lim_{t\uparrow \zeta^y}X^{D,y}_t =y \qquad\P_x^y\mbox{-a.s.}.
\end{equation}

Let $\{D_k, \, k \ge 1\}$ be an increasing sequence of
relatively compact open subsets of $D$ such that
$D_k \subset \overline{D_k} \subset D$ and
$\cup_{k=1}^\infty D_k =D$. Then
$$
\E_x \left[
\overline{G}_D(X^D_{\tau_{D_k}}, y) \right]
=  \overline{G}_D(x, y)  \P_x^y ( \tau_{D_k}<\zeta^y).
$$
By (\ref{e_lim}),
we have $\lim_{k \to \infty} \P_x^y ( \tau_{D_k}<\zeta^y) =0$. Thus
\begin{equation}\label{eqn:5.4}
\lim_{k \to \infty}  \E_x \left[
\overline{G}_D(X^D_{\tau_{D_k}}, y) \right] =0.
\end{equation}
The last claim in the theorem is easy. By (\ref{est:6.1}) and
(\ref{G_bd}),
for every $x,y \in D$, we have
$$
\E_x\left[\,\overline{G}_D(X^D_t, y)\right] \le
\frac{c}{t^{\frac{d}{2}}h_D(y)}
\int_D \frac{dz}{|z-y|^{d-2}},
$$
which converges to zero as $t$ goes to $\infty$.
\qed

\medskip

The proof of the next proposition can be found in the proofs of
Theorems 2-3 in \cite{P}. Since we will refer to the argument of the proof of
the proposition later,
we include  the proof here for the reader's convenience.

\medskip

\begin{prop}\label{p:ph}
If $h$ is a nonnegative harmonic function for $X^D$ and $U$ is
an open subset of $D$
with $\overline{U} \subset D$,
then there exists a Radon measure $\nu$ supported on $\partial U$
such that $h= G_D \nu$ in $U$. In particular, every
nonnegative harmonic function
for $X^D$ is continuous.
\end{prop}

\pf
Using (\ref{G_bd}) and (A1)-(A2), we see that the properties (a)-(b) in
\cite{P2} is true for $X^D$. Thus by Corollary 1 in \cite{P2},
$h$ is excessive.
Let $T_U:=\inf \{ t > 0: X^D_t \in U\}$.
Since $h$ is excessive,
Corollary 1 to Theorem 2 in \cite{CR} implies that there exists
a Radon measure $\nu$ supported on $\overline{U}$ such that
$\E_x[h(X_{T_U}^D)]= G_D \nu(x)$ for all $x \in D$. Since
$$
{G}_D \nu (x)\,=\,
\int_{U}\overline{G}_D (x,y) \nu(dy) +\int_{  \partial{U}  }
\overline{G}_D (x,y) \nu(dy) =: h_1(x)+h_2(x),
\qquad x \in D
$$ and  $h_1$ and $h_2$ are excessive (Theorem \ref{t:Riesz}),
$h_1$ and $h_2$ must be harmonic with respect to $X^D$ .
Let $K$ be a compact subset of $U$. By the harmonicity of $h_1$, we have
$$
\E_x[h_1(X_{T_{K^c}}^D)]= \int_{  \partial{U}  }
\overline{G}_D (x,y) \nu(dy).$$
But, by Corollary 1 to Theorem 2 in \cite{CR}, $\nu$ can not charge
the interior of $K$.
Since $K$ is an arbitrary compact subset of $U$,
we get that $h_1$ is identically zero and $\nu$ is supported by
$\partial U$. Therefore
we have shown $h(x)=\E_x[h(X_{T_U}^D)]=G_D \nu(x)$ for $x \in U$.
Now the continuity of $h$ follows from the continuity of
$\overline{G}_D (x,y)$.
\qed

\section{Dual of $X^D$ in bounded domains}

In this section we assume that
$D$ is an arbitrary bounded domain.
First we show that $X^D$ has a nice dual process with respect to $\xi_D$
and then we will study
the dual process of $X^D$.

Recall that
$
h_D(x)=\ \int_D G_D(y,x) dy$, $ \xi_D(dx)=h_D(x) dx$,
$\overline{q}^D(t,x,y)= \frac{ q^D(t,x,y)}{h_D(y)}$ and
$ \overline{G}_D(x,y)= \frac{ G_D(x,y)}{h_D(y)}$.
We note that
$$
\int_D \overline{G}_D (x,y) \xi_D (dx) \le
\frac{\|h_D\|_{L_{\infty}(D)}}{h_D(y)}
\int_D G_D (x,y) dx = \|h_D\|_{L_{\infty}(D)} < \infty.
$$
So we have
\begin{description}
\item{(A5)} for every compact subset $K$ of $D$, $\int_K
\overline{G}_D(x,y) \xi_D(dx) < \infty$.
\end{description}
(A1)-(A5), (\ref{G_bd}) and Theorem \ref{t:p} imply that
the conditions (i)-(vii) and (70)-(71) in
\cite{L1} (also see Remark on page 391
in \cite{L}) are satisfied,
thus $X^D$ has  a continuous Hunt process as a dual process.

\medskip

\begin{thm} \label{dualpr}
There exists a transient continuous Hunt process $\wh X^D $
in $D$ such that $\wh X^D $
is a strong dual of $X^D$ with respect to the measure $\xi_D$,
that is, the density of the semigroup
$\{ \wh P^D_t\}_{t\geq 0}$ of $\wh X^D$ is given by
$\wh q^D (t,x,y) := \overline{q}^D(t,y,x)$ and thus
$$
\int_D f(x) P^D_t g(x) \xi_D(dx)=\int_D g(x) \wh P^D_t f(x) \xi_D(dx)
\quad \hbox{\rm for all } f, g\in L^2(D, \xi_D).
$$
\end{thm}

\pf The existence of a dual continuous Hunt process $\wh X^D$ is proved in
\cite{L1}. To show $\wh X^D$ is transient, we need to show that for every
compact subset  $K$ of $D$, $\int_K \overline{G}_D (x,y) \xi_D (dx)$
is bounded. This is just (A5) above.
\qed

\medskip

We will use $\wh{\zeta}$ to denote the lifetime of $\wh X^D$. Note
that the process $\wh X^D$ might have killing inside $D$, that is,
the probablity of the event $\wh X^D_{{\wh \zeta}-}\in D$ might be
strictly positive.

By (A1), (A2) and (A5), $\wh X^D$ is a transient diffusion
satisfying the conditions in \cite{P} and \cite{CR}.
So by applying  Theorem 1 in \cite{P}, we have that
 \begin{description}
\item{(A6)}
 for each $y$,  $ x \to \overline{G}_D(y,x)$ is excessive for
$\wh X^D$ and harmonic in
 $D \setminus \{y\}$.
Moreover, for every
open subset $U$ of $D$, we have
\begin{equation}\label{e:GH0}
\E_x[\overline{G}_D(y, \wh X^D_{\wh T_{U}})]=\overline{G}
_D(y,x), \qquad (x,y) \in D \times U
\end{equation}
 where  $\wh T_U:=\inf \{ t > 0: \wh X^D_t \in U\}$. In particular,
for every $y \in D$ and $\eps > 0$,
 $\overline{G}_D( y, \,\cdot\,)$
is regular harmonic in $D \setminus B(y, \eps)$
with respect to $\wh X^D$.
\end{description}

Thus the Riesz representation theorem (Theorem \ref{t:Riesz}) is
valid for $\wh X^D$ too.

\medskip

\begin{thm}\label{t:Riesz_dual}
\begin{description}
\item{(1)}  If $u$ is a potential for $\wh X^D$, then
there exists
a unique Radon measure $\nu$ on $D$ such that
$$
u(x) \,=\,\wh G_D \nu (x)\, := \int_D\overline{G}
_D(y,x) \nu (dy) $$
\item{(2)} If $f$ is an excessive function for $ \wh X^D$
and $f$ is not identically zero,
then there exists
a unique Radon measure $\nu$ on $D$ and a nonnegative harmonic function $h$
for $\wh X^D$ such that
$ f=\wh G_D \nu + h.$
\end{description}
\end{thm}

\medskip

\begin{thm}\label{t:p_dual}
 For each $y$,  $ x \to \overline{G}_D(y,x)$ is a pure potential for
$ \wh X^D$. In fact,
for every sequence $\{U_n\}_{n \ge 1}$
of open sets with
$\overline{U_n} \subset U_{n+1}$ and  $\cup_n U_n = D$,
$$
 \lim_{n \to \infty}
 \E_x\left[\overline{G}_D(y, \wh X^D_{\wh\tau_{U_n}});
\wh\tau_{U_n}<\wh\zeta\right]\,=\, 0,
\qquad x \not= y.
$$
Moreover, for every $x,y \in D$, we have
$$
\lim_{t \to \infty} \E_x\left[\overline{G}_D(y, \wh X^D_t);
t<{\wh\zeta}\right] \,=\, 0.
$$
\end{thm}
\pf
The first assertion can be proved using an argument same as that of
the proof of
Theorem \ref{t:p}, so we only need to prove the last assertion.

By (\ref{est:6.1}) and
(\ref{G_bd}), we have for every $x,y \in D$,
$$
\E_x\left[\overline{G}_D(y,  \wh X^D_t);
t<{\wh\zeta}\right]
\,=\,
\int_D \frac{q^D(t,z,x)}{h_D(x)} \overline{G}_D(y,z) h_D(z)dz
\,\le\, \frac{c}{t^{\frac{d}{2}}h_D(x)}
\int_D \frac{dz}{|z-y|^{d-2}},
$$
which converges to zero as $t$ goes to $\infty$.
\qed

\medskip

Note that  every nonnegative harmonic function for $ \wh X^D$ is
excessive and continuous by Corollary 1
 in \cite{P2}.
Let  $\{\wh G^D_\alpha, \, \alpha \geq 0\}$ be the resolvent of
$ \wh X^D$ with respect $\xi_D$ so that
$\wh G^D_0 (y, dx)=\overline{G}_D (x, y) \xi_D(dx)$.

\medskip

\begin{prop}\label{s_feller_kd}
$\wh X^D$  has the strong Feller property in the resolvent sense;
that is,
for every bounded Borel function $f$ on $D$ and $\alpha \ge 0$,
$\wh G^D_\alpha f(x)$
is bounded continuous function on $D$.
\end{prop}

\pf
By the resolvent equation $\wh G^D_0 =\wh G^D_\alpha +\alpha
\wh G^D_0 \wh G^D_\alpha$,
it is enough to show the strong Feller property for $\wh G^D_0$.
Fix  a bounded Borel  function $f$ on $D$ and a sequence
$\{ y_n \}_{n \ge 1}$ converges to $y$ in $D$.
Let $M := \| f h_D \|_{L_\infty(D)} <\infty$.
We assume  $\{ y_n \}_{n \ge 1} \subset K$ for a compact subset $K$ of $D$.
Let $A:= \inf_{y \in K} h_D (y)$.
By Proposition \ref{e_m}, we know that
$A$ is strictly positive.
Note that there exists a constant $c_1$ such that for every $\delta>0$
$$
\left(\int_{B(y,\delta)} \frac{dx}{|x- y|^{d-2}}
+\int_{B(y_n,2\delta)} \frac{dx}{|x- y_n|^{d-2}} \right) \le c_1 \delta^2.
$$
Thus by (\ref{G_bd}), there exists a constant $c_2$ such that for
every $\delta>0$ and
$y_n$ with $y_n \in  B(y, \frac{\delta}2) \subset  B(y, 2\delta) \in K$,
\begin{eqnarray*}
&&
\int_{B(y,\delta)} \overline{G}_D (x, y)  f(x) \xi_D(dx)
+\int_{B(y,\delta)} \overline{G}_D (x, y_n)  f(x) \xi_D(dx)\\
&& \le \frac{M}{A}\left(\int_{B(y,\delta)}G_D (x, y) dx
+\int_{B(y_n,2\delta)} G_D (x, y_n) dx\right)\\
&& \le \frac{c_2M}{A}\left(\int_{B(y,\delta)} \frac{dx}{|x- y|^{d-2}}
+\int_{B(y_n,2\delta)} \frac{dx}{|x- y_n|^{d-2}} \right)\, \le\,
\frac1{A}c_1c_2M \delta^2
\end{eqnarray*}
Given $\eps$, choose $\delta$ small enough such that
$\frac1{A}c_1c_2 M \delta^2 < \frac{\eps}2$. Then
$$
|\wh G^D_0 f(y)-\wh G^D_0 f(y_n)|\, \le\, M  \int_{ D \setminus B(y,\delta)}
| \overline{G}_D (x, y) - \overline{G}_D (x, y_n) | dx \,+\,   \frac{\eps}2.
$$
Note that  $\overline{G}_D (x, y_n)$ converges to
$\overline{G}_D (x, y)$ for every $x\not=y$
and $\{\overline{G}_D (x, y_n)\}$ are uniformly
bounded on $x \in  D \setminus
B(y,\delta)$ and $y_n \in B(y, \frac{\delta}2)$.
So the first term on the right hand side of the inequality
above goes to zero as $n \to \infty$ by the bounded convergence
 theorem.
\qed

\section{Comparison of harmonic measures and scale invariant
Harnack inequality for the dual process}

In this section we still assume that $D$ is an  arbitrary bounded domain.
For any open subset $U$ of  $D$,
we use $\wh X^{D, U}$ to denote
the process obtained by killing $\wh X^D$ upon exiting $U$,
i.e., $\wh X^{D, U}_t(\omega)=\wh X^D_t(\omega)$ if $t< \wh\tau^D_U(\omega)$
and $\wh X^{D,U}_t(\omega)=\partial$ if $t\geq \wh \tau^D_U(\omega)$,
where $\wh \tau^D_U:=\inf\{t >0 : \wh X^D_t \notin U\}$ and $\partial$
is the cemetery state. Then by Theorem 2 and Remark 2 after it in \cite{S},
$X^U$
and $\wh X^{D, U}$ are dual processes with respect to $\xi_D$.
Now we let
$$
\wh q^{D,U}(t,x,y)\,:=\,\frac{q^U(t,y,x)h_D(y)}{h_D(x)}.
$$
By the joint
continuity of $q^D(t,x,y)$ (Theorem 2.4 in \cite{KS}) and the
continuity and positivity of $h_D$ (Proposition \ref{e_m}), we know that $q^{D,U}(t,\cdot,\cdot)$
is jointly continuous on $U \times U$. Thus we have the following.

\medskip

\begin{thm}\label{t:kill_dual}
For every open subset $U$, $\wh q^{D,U}(t,x,y)$
is jointly continuous on $U \times U$ and is the transition density of $\wh X^{D, U}$ with respect to
the Lebesgue measure. Moreover,
\begin{equation}\label{e:dGk}
\wh G_{D,U}(x,y)\,:=\,\frac{G_U(y,x)h_D(y)}{h_D(x)}
\end{equation}
is the Green function of $\wh X^{D, U}$ with respect to the Lebesgue measure
so that for every nonnegative Borel function $f$,
$$
\E_x\left[\int_0^{\wh \tau^D_U} f\left(\wh X^{D}_t\right)dt\right]\,=\,
\int_U \wh G_{D,U}(x,y) f(y) dy.
$$
\end{thm}

\medskip

Using (\ref{e:dGk}), one can check that $\wh X^{D, U}$ satisfies the
conditions
in \cite{CR} and \cite{P}. Thus by repeating the argument in the
proof of Proposition \ref{p:ph}, we get the following.

\medskip

\begin{prop}\label{p:phk}
If $h$ is a nonnegative harmonic  for $\wh X^D$ in $U$
and $V$ is an open subset of $U$
with $\overline{V} \subset D$,
then there exists a Radon measure $\nu$ supported on $\partial V$
such that
$$
h(x)= \int_{\partial V} \frac{G_U(y,x)h_D(y)}{h_D(x)}
\nu(dy), \quad  x \in V.
$$
In particular, every nonnegative harmonic function
for $\wh X^D$ in $U$ is continuous.
\end{prop}

\medskip

Using (\ref{e:2sge}) and Proposition \ref{e_m}, we see that
for every compact subset $K$ of $D$, there exist $c_1$, $c_2$
and $c_3$ such that
 for every positive $t_0$ and $\delta$,
\begin{eqnarray}
\sup_{t\le t_0,x \in K, |x-y| >  \delta} \frac{q^D(t , y,x) h_D(y)}{h_D(x)}
&\le& c_1 e^{c_2 t_0}
\sup_{t \le t_0,x \in K, |x-y| >  \delta}  t^{-\frac{d}{2}}
e^{-c_3\frac{|x-y|^2}{t}} \nonumber\\
&\le& c_1 e^{c_2 t_0}
\sup_{t \le t_0}  t^{-\frac{d}{2}}
e^{-c_3\frac{\delta^2}{t}} \,<\, \infty.\label{e:bd1}
\end{eqnarray}
(\ref{e:bd1}) implies that for any compact subset $K$ of $D$,
\begin{eqnarray*}
&&\sup_{t\le t_0, x \in K} \P_x(| \wh X^D_t -x| \ge \delta;
t<\wh\zeta)
\le   c_1 e^{c_2 t_0} \sup_{t\le t_0, x \in K} \int_{ |x-y| \ge \delta}
t^{-\frac{d}{2}}
e^{-c_3\frac{|x-y|^2}{t}} dy\\
&&=  c_4 e^{c_2 t_0}  \sup_{t\le t_0} \int^\infty_{\delta}
t^{-\frac{d}{2}} r^{d-1}
e^{-c_3\frac{r^2}{t}} dr
\le c_5  e^{c_2 t_0}  \int^\infty_{\frac{\delta} {\sqrt{t_0}}}
u^{d-1} e^{-c_3 u^2} du
\end{eqnarray*}
for some $c_5=c_5(d)$ and $c_6=c_6(d)$ .
Thus
\begin{equation}\label{e:con}
\lim_{t_0 \downarrow 0}  \sup_{t\le t_0, x \in K}
\P_x(|\wh X^D_t -x| \ge \delta; t<\wh\zeta) \,=\,
\lim_{t_0 \downarrow 0}  \sup_{t\le t_0, x \in K}
\P_x(\wh X^D_t \in D \setminus B(x, \delta)) \,=\,0.
\end{equation}

Using (\ref{e:con}) we can easily prove the next lemma.

\medskip

\begin{lemma}\label{l:con}
For any $\delta>0$ and $x \in D$ with $B(x, 2\delta) \in D$, we have
$$
\lim_{s\downarrow 0} \sup_{x \in D:B(x, 2\delta) \in D}
 \P_x(\wh  \tau^D_{B(x, \delta)} \le s\wedge\wh\zeta) \, =0.
$$
\end{lemma}

\pf For any $t>0$ and any Borel set $A$ in $D$, we put
$N_t(x, A)=\P_x(\wh X^D_t\in A)$.
Then by an extended version of the strong Markov property
(see page 43--44 of \cite{BG}), we have for every $ x \in D $ with
$B(x, 2\delta) \in D$,
\begin{eqnarray*}
&& \P_x( \wh \tau^D_{B(x, \delta)} \le s\wedge\wh\zeta) \le
\P_x\left( \wh  \tau_{B(x, \delta)} \le s, \,\wh X_s \in B(x, \frac{\delta}2)
\right)\,+\, \P_x\left(
\, \wh X^D_s \in B(x, \frac{\delta}2)^c, s<\wh\zeta \right)\\
&&\le  \E_x\left[N_{s-\wh \tau^D_{B(x, \delta)}}\left(\wh
X^D_{\wh\tau_{B(x, \delta)}},
B(x, \frac{\delta}2)\right);\,
\wh \tau^D_{B(x, \delta)} \le s \right]
 \,+\, \P_x\left(    \,\wh X^D_s \in B(x, \frac{\delta}2)^c,
s<\wh\zeta \right)
\end{eqnarray*}
Since $X^D_{\tau_{B(x, \delta)}} \in \partial B(x, {\delta})$
almost surely on $\{\wh\tau_{B(x, \delta)}<\wh\zeta\}$,
the conclusion of the lemma follows from (\ref{e:con}).
\qed

\medskip
A  bounded domain $U$ in $\R^d$ is said to be
a $C^{1,1}$
domain if there is a localization radius
$r_0>0$
 and a constant
$\Lambda >0$
such that for every $Q\in \partial U$, there is a
$C^{1,1}$-function $\phi=\phi_Q: \R^{d-1}\to \R$ satisfying $\phi (0)
= \nabla\phi (0)=0$,
$\| \nabla \phi  \|_\infty \leq \Lambda$,
$| \nabla \phi (x)-\nabla \phi (z)| \leq \Lambda
|x-z|$, and an orthonormal coordinate
system $y=(y_1, \cdots, y_{d-1}, y_d):=(\tilde y, y_d)$
such that
$ B(Q, r_0)\cap D=B(Q, r_0)\cap \{ y: y_d > \phi (\tilde y) \}$.

Using (\ref{e:dekbmsd}), it is easy to show the following.

\medskip

\begin{prop}\label{s_feller}
For any bounded
$C^{1,1}$ domain $U\subset D$ with $\overline{U} \subset D$,
$\wh X^{D,U}$ satisfies the strong Feller property in the semigroup sense;
that is,
for every bounded Borel function $f$ on $U$ ,
$$
\E_x\left[f(\wh X^D_t);t < \wh \tau^D_{U}\right]
$$
is a bounded continuous function on $U$.
\end{prop}

\pf
Fix $x_0 \in U$ and $t >0$.
Suppose $x_n \in U$ converges to $x_0 \in U$. Let
$N:=\inf_{n \ge 1} h_D(x_n) > 0$
Then
\begin{eqnarray*}
&&\left|\E_{x_n}\left[f(\wh X^D_t);t < \wh \tau^D_{U}\right]
- \E_{x_0}\left[f(\wh X^D_t);t < \wh \tau^D_{U}\right]\right|\\
&\le&c_1 \int_U
\left| \frac{q^U(t,y,x_n)}{h_D(x_n)} -\frac{q^U(t,y,x_0)}{h_D(x_0)}   \right|
h_D(y)dy\\
&\le&c_2 \int_U \left(
\left| \frac{q^U(t,y,x_n)}{h_D(x_n)}-\frac{q^U(t,y,x_0)}{h_D(x_n)}\right|+
q^U(t,y,x_0)\left| \frac1{h_D(x_n)}-\frac1{h_D(x_0)}   \right| \right)
dy\\
&\le&\frac{c_2}{N} \int_U
\left| q^U(t,y,x_n)-q^U(t,y,x_0)\right| dy+ c_2
\left|\frac1{h_D(x_n)}
-\frac1{h_D(x_0)}   \right|.
\end{eqnarray*}
Given $\eps>0$, choose $n_0 >0$ such that
$$
c_2\left|\frac1{h_D(x_n)}
-\frac1{h_D(x_0)}   \right| < \frac{\eps}{4}, \quad n \ge n_0.
$$
Let $\rho_U(y)$ be the distance between $y$ and $\partial D$.
By (\ref{e:dekbmsd}),
\begin{eqnarray*}
&&\int_{U } \left| q^U(t,y,x_n)-q^U(t,y,x_0)\right| dy\\
&&=\,\left(\int_{ \{y \in U: \rho_U(y)< \delta\}}+
\int_{   \{y \in U: \rho_U(y) \ge \delta\} }\right)
\left| q^U(t,y,x_n)-q^U(t,y,x_0)\right| dy\\
&&\le\, c_3\,|U| \,\delta t^{-\frac{d+1}{2}}\,+\,
\int_{ \{y \in U: \rho_U(y) \ge \delta\} }
\left| q^U(t,y,x_n)-q^U(t,y,x_0)\right| dy,
\end{eqnarray*}
for some $c_3$.
Now we choose $\delta$ small so that
$c_2 c_3|U| N^{-1} \delta t^{-\frac{d+1}{2}} <\frac{\eps}4. $
The convergence of the second term on the right hand side of the inequality
above follows
from the uniform continuity of
$q^{U}(t, \,\cdot\, , \,\cdot \,)$ on $B(x_0, \frac12\rho_U(x_0))
\times  \{y \in U: \rho_U(y) \ge \delta\}$
(Theorem 3.1 in \cite{KS}).
Thus we have proved the proposition.
\qed

\medskip

Recall that, a point $z$ on the boundary $\partial U$ of an open
subset $U$ of $D$
is said to be a regular boundary point for $\wh X^D$ in $U$
if $\P_x( \wh \tau^D_U=0)=1$.
An open subset $U$ of $D$ is said to be regular
if every point in $\partial U$
is a regular boundary point.

\medskip

\begin{prop}\label{p:reg}
Suppose $U$ is an open subset of $D$ with $\overline{U} \subset D$,
and $z\in \partial U$.
If there is a cone $A$ with vertex $z$ such that
$A\cap B(z, r)\subset D \setminus U$ for some $r>0$, then $z$ is a regular
boundary point of $U$.
\end{prop}

\pf
Choose a bounded smooth domain $D_1$ with $\overline{U} \subset D_1
\subset \overline{D_1} \subset D$.
Without loss of generality, we may assume that $z=0$ and
$A\cap B(z, r)\subset D_1 \setminus U$.
For $n\ge 1$, put $r_n=r/n$. Under $\P_0$, we have
$$
\bigcap^{\infty}_{m=1}\bigcup^{\infty}_{n=m}
\{\wh X^D_{r_n}\in A\cap B(0, r)\}
\subset \{\wh \tau^D_U=0\}.
$$
Moreover, since $D_1$ is bounded smooth and $\overline{D_1} \subset D$,
 by (\ref{e:dekbmsd}) there exists a constant $c_1>0$ such that
for $x \in A\cap B(0, r)$ and large $n$
$$
\frac{q^{D_1}(t,x,0)h_D(x)}
{h_D(0)} \ge c_1 r_n^{-\frac{d}2}e^{-\frac{c_2|x|^2}{r_n}}.
$$
Hence
\begin{eqnarray*}
&&\P_0(\wh \tau^D_U=0)
\,\ge\,\P_0\left(\bigcap^{\infty}_{m=1}\bigcup^{\infty}_{n=m}
\{\wh X^D_{r_n}\in A\cap B(0, r)\}\right)\,
\ge\,\limsup_{n\to\infty}\P_0(\wh X^D_{r_n}\in A\cap B(0, r))\\
&&\ge\limsup_{n\to\infty}\P_0(\wh X^{D_1}_{r_n}\in A\cap B(0, r))
\,=\, \limsup_{n\to\infty} \int_{A\cap B(0, r)} \frac{q^{D_1}(t,x,0)h_D(x)}
{h_D(0)}dx\\
&&\ge\, \limsup_{n\to\infty}
c_1
\int_{A\cap B(0, r)}r_n^{-\frac{d}2}e^{-\frac{c_2|x|^2}{r_n}}dx
\,\ge\, \limsup_{n\to\infty}
c_1\int_{A\cap B(0, n)}e^{-c_2|y|^2}dy\,>\,0.
\end{eqnarray*}
The assertion of the proposition now follows from Blumenthal's
zero-one law (Proposition I.5.17 in \cite{BG}).
\qed

\medskip

This result implies that all bounded Lipschitz domains
(see below for the definition) are regular if their closures are in $D$.
Modifying the argument in the second part of the proof of Theorem 1.23
in \cite{CZ}, we get the following result.

\medskip

\begin{prop}\label{p:cont_b}
Suppose $U$ is an open set subset of $D$
with $\overline{U} \subset D$ and $f$ is a bounded
Borel function on $\partial U$. If $z$ is a regular boundary point of $U$ for
$\wh X^D$ and $f$ is continuous at $z$
$$
\lim_{\overline{U}\ni x\to z}\E_x\left[f\left(\wh X^D_{\wh \tau^D_U}
\right);
\wh \tau^D_U   <\wh\zeta\right]=f(z).
$$
\end{prop}

\pf
With Lemma \ref{l:con} and Proposition \ref{s_feller} in hand,
the proof is routine. We omit the details.
\qed

\medskip

Let $G^{0}_D$ be the Green function of a
Brownian motion $W$ in $D$.
By Theorem 3.7 in \cite{KS1},  there exist constants $r_1=r_1(d, \mu)>0$ and
$M_1=M_1(d, \mu) >1$ depending on
$\mu$ only via the rate at which $\max_{1 \le i \le d}M_{\mu^i}(r)$
goes to zero such that
for $r \le r_1$, $z \in \R^d$, $x, y\in B(z,r)$,
\begin{equation}\label{e:Green_B}
M_1^{-1}\, G^0_{B(z,r)}(x,y) \le\, G_{B(z,r)}(y,x)\,
\le\, M_1\, G^0_{B(z,r)}(x,y).
\end{equation}
We will fix the constants $r_1>0$ and $M_1>0$ above in the
remainder of this section.

\medskip

\begin{thm}\label{t:7.5}
For any bounded domain $D$, $r \le r_1 $, $z \in D$ and
$x \in B^z_r:=B(z,r) \in \overline{B(z,r)}
\subset D  $, we have
\begin{equation}\label{h_est}
M_1^{-2}
 \,h_D(y)\P_x\left( W_{\tau_{B^z_r}} \in dy\right)
\,\le   \,h_D(x)\P_x\left( \wh X^D_{\wh\tau_{B^z_r}} \in dy,
\wh\tau_{B^z_r}<\wh\zeta\right)
\,\le \, M_1^2 \,h_D(y)\P_x\left( W_{\tau_{B^z_r}} \in dy\right).
\end{equation}
\end{thm}

\pf We fix $z_0 \in D $ and $r \le r_1 $ with
$ B(z_0,r) \in \overline{B(z_0,r)} \subset D  $, and let $B:=B(z_0,r)$.
The idea of the proof is similar to that of Theorem 2.2 in \cite{CWZ}. We
include the detail here for the reader's convenience.
Let $\varphi \ge 0$ is a continuous function on $\partial B$ and
let
$$
u(x):= \E_{x} [ \varphi( \wh X^D_{\wh \tau^D_B});\wh \tau^D_B<\wh\zeta ].
$$
By Proposition \ref{p:cont_b}, we know that
$u$ is harmonic for $\wh X^D$ in $B$ and continuous on $\overline{B}$.
Let $B(n):=B(z_0,(1-\frac1{n})r)$,
$T_{n}:=\inf \{ t > 0: \wh X^D_t \in B(n)\}$ and
$
u_{n}(x):= \E_{x}  [ u(\wh X^{D,B}_{T_n})].
$
Then by Proposition \ref{p:phk}, there exist Radon measures $\nu_n$
supported on $\partial B(n) $  such that
$$
u_{n}(x) \,=\,\frac1{h_D(x)}\int_{\partial B(n)}
G_B(y,x)  \nu_n (dy).
$$
Let
$$
v_{n}(x) \,:=\,\int_{\partial B(n)} G^{0}_B (x,y)
\nu_n (dy).
$$
Then by (\ref{e:Green_B}),
$$
M_1^{-1} v_n (x) \,\le\,  h_D(x) u_n(x) \,\le\, M_1 v_n (x),
\qquad x \in B(n).
$$
Since $v_n$ is a harmonic function in $B(m)$ for $m \ge n$, by the H\"{o}lder
continuity of $v_n$ and a diagonalization procedure, there is
a subsequence $n_k$ such
that $v_{n_k}$ converges uniformly on each $B(m)$ to a harmonic
function $v$ in $B$. Thus
\begin{equation}\label{e:co}
M_1^{-1} v (x) \,\le\,  h_D(x) u(x) \,\le\, M_1 v (x), \qquad x \in B.
\end{equation}
Since $B$ is regular for $\wh X^D$ (Proposition \ref{p:reg}), by taking the
limit above and using Proposition \ref{p:cont_b},
we get for every $w \in \partial B$
\begin{equation}\label{e:coo}
M_1^{-1} h_D(w) \varphi (w)  \,\le\, \liminf_{B \ni x \to w}  v (x)
\,\le\, \limsup_{B \ni x \to w}  v (x)
\,\le\, M_1 h_D(w) \varphi (w).
\end{equation}
Let
$$
w(x)= \E_x\left[h_D( W_{\tau_{B}} ) \varphi (W_{\tau_{B}}) \right].
$$
$w$ is a harmonic function in $B$ and continuous in $\overline{B}$ with
the boundary value $ h_D(w) \varphi (w)$.
Thus by the maximum principle and (\ref{e:coo}), we get
$M_1^{-1} w(x) \le v (x) \le M_1 w(x) $. So by (\ref{e:co})
$M_1^{-2} w(x)\le h_D(x) u(x) \le M_1^2 w(x)$, which is
\begin{eqnarray*}
M_1^{-2} \int_{\partial B} \varphi (w) h_D(w)  \P_x\left(
W_{\tau_{B}} \in dw\right)
&\le& \int_{\partial B} \varphi (w) h_D(x)   \P_x\left(
\wh X^D_{\wh \tau^D_{B}} \in dw, \wh \tau^D_{B}<\wh\zeta \right)  \\
&\le&M_1^{2} \int_{\partial B} \varphi (w) h_D(w)  \P_x\left(
W_{\tau_{B}} \in dw\right).
\end{eqnarray*}
\qed

\medskip

Let $\rho_D(x)$ be the distance between $x$ and $\partial D$.

\medskip

\begin{lemma}\label{l:rho}
Suppose $D$ is a bounded $C^{1,1}$ domain. Then there exists a
constant $c=c(D)$ such that
\begin{equation}\label{e:rho}
\frac1{c}  \rho_D(x) \,\le \, h_D(x) \,\le \, c \rho_D(x).
\end{equation}
\end{lemma}

\pf
Since $D$ is bounded, the Green function estimates for
$X^D$ ((6.2)-(6.3) in \cite{KS}) imply that
$$
c_1 \rho_D(x) (1\wedge \frac{\rho_D(y)} {|x-y|^2})
\frac{1}{|x-y|^{d-2}}
\,\le
\, G_{D}(y,x) \,\le\, c_2\,
 \frac{\rho_D(x)}{|x-y|^{d-1}}
$$
for some positive constants $c_1$ and $c_2$.
Integrating over $y$ we get
$$
c_1 \rho_D(x) q_1(x) \,\le \, h_D(x) \,\le \, c_2\rho_D(x) q_2(x)
$$
where
$$
q_1(x):= \int_D (1\wedge \frac{\rho_D(y)} {|x-y|^2})
\frac{1}{|x-y|^{d-2}}dy \quad \mbox{and} \quad
q_2(x):=\int_D     \frac{1}{|x-y|^{d-1}}            dy.
$$
By elementary calculus, we easily see that
$\inf_{x \in D} q_1(x) >0$ and $\sup_{x \in D} q_2(x) < \infty$.
\qed

\medskip

Let $\psi_0 $ be the ground state for the killed Brownian motion in $D$,
that is, $\psi_0$ is the positive eigenfunction corresponding
to the smallest eigenvalue of the Dirichlet Laplacian $-\frac12\Delta|_D$
with $\int_D\psi^2_0(dx)dx=1$.
If $D$ is bounded $C^{1,1}$,
it is well-known that there exists $c_1$ such that $ c_1^{-1}
\rho_D(x)\le \psi_0(x) \le
c_1\rho_D(x)$. So we get the next result as a corollary of
 Theorem \ref{t:7.5} and Lemma \ref{l:rho}.

\medskip

\begin{corollary}\label{g_com}
For any bounded $C^{1,1}$ domain $D$,
there exists a positive constant $c$ such that for every
$r \le r_1 $, $z \in D$ and $x \in B^z_r:=B(z,r) \in \overline{B(z,r)}
\subset D  $, we have
\begin{equation}\label{h_est1}
c^{-1}
 \,\P_x\left( W^{\psi_0}_{\tau_{B^z_r}} \in dy\right)
\,\le  \P_x\left( \wh X^D_{\wh\tau_{B^z_r}} \in dy,
\wh\tau_{B^z_r}<\wh\zeta\right)
\,\le \, c\,\P_x\left( W^{\psi_0}_{\tau_{B^z_r}} \in dy\right)
\end{equation}
where $W^{\psi_0}$ is  the $h$-conditioned Brownian motion
with $h=\psi_0$.
\end{corollary}

\medskip

In the remainder of this section, we will prove a scale
invariant Harnack inequality for
$\wh X^D$ for any bounded domain $D$. First we prove the
following lemma. Recall
that $r_1>0$ and $M_1>0$ are the constants from (\ref{e:Green_B}).

\medskip

\begin{lemma}\label{l:dis}
There exists  a constant $c=c(D,\mu) >1$ such that
for every $r < r_1$ and $B(z, r)$ with $B(z, r)
\subset D$,
$$
\frac{h_D(x)}{h_D(y)} \le c, \quad x,y \in B(z, \frac{r}2).
$$
\end{lemma}

\pf Since $r < r_1$, by (\ref{G_bd}) and (\ref{e:Green_B}),
there exists $c_1=c_1(D, \mu) >1$ such that
for every $x, w \in \overline{B(z, \frac{3r}4)}$
$$
c_1^{-1}\, \frac{1}{|w-x|^{d-2}} \,\le\, M_1^{-1}\, G^{0}_{B(z,r)}
(w,x) \,\le\,  G_{B(z,r)}  (w,x)
\,\le\,  G_{D}  (w,x) \,\le\, c_1 \,\frac{1}{|w-x|^{d-2}}.
$$
Thus for $ w \in \partial {B(z, \frac{3r}4)}$ and $ x,y \in
{B(z, \frac{r}2)}$, we have
\begin{equation}\label{e:pp1}
G_{D}  (w,x) \,\le\, c_1 \,\left(\frac{|w-y|}{|w-x|}\right)^{d-2}
\frac{1}{|w-y|^{d-2}} \,\le\, 4^{d-2}\, c_1^2\, G_{D}  (w,y).
\end{equation}
On the other hand, from (\ref{e:GH}), we have
\begin{equation}\label{e:pp2}
h_D(x) \,=\, \int_D G_D(a,x) da\, = \,
\int_D \E_a\left[G_D(X_{T_{B(z,\frac{3r}4)}},x)\right] da ,
\quad x \in B(z,\frac{3r}4).
\end{equation}
Since $X_{T_{B(z,\frac{3r}4)}} \in \partial {B(z, \frac{3r}4)}$,
combining (\ref{e:pp1})-(\ref{e:pp2}), we get
$$
h_D(x) \,\le\,  4^{d-2}\, c_1^2\,
\int_D \E_a\left[G_D(X_{T_{B(z,\frac{3r}4)}},y)\right] da
\,=\, 4^{d-2}\, c_1^2 \,h_D(y), \quad x, y \in B(z,\frac{r}2)
$$
\qed

\medskip

Now we are ready to prove the scale invariant Harnack inequality.

\medskip

\begin{thm}\label{t:HP} (Scale invariant Harnack inequality)
There exists $N=N(d, \mu) >0$ such that for  every
harmonic function $f$ of $\wh X^D$ in $B(z_0, r)$ with
$r \in (0 , r_1]$ and $B(z_0, r) \subset D$,  we have
$$
\sup_{y \in B(z_0, r/4)}f(y) \le N \inf _{y \in B(z_0, r/4)}f(y)
$$
\end{thm}

\pf We fix  $z_0 \in D $ and $r \le r_1 $ with
$ B(z_0,r) \in B(z_0,r) \subset D  $, and a harmonic function $f$
of $\wh X^D$ in $B(z_0, r)$.
By the harmonicity of $f$ and Theorem \ref{t:7.5}, for every
$x \in B(z_0,\frac{r}2)$
$$
f(x)\,=\, \E_x\left[f\left(\wh X^D_{\wh\tau_{B(z_0,\frac{r}2)}}\right);
\wh\tau_{B(z_0,\frac{r}2)}<\wh\zeta \right]
\,\le\, \frac{M_1^2}{h_D(x)}\E_x\left[h_D\left(
W_{\tau_{B(z_0,\frac{r}2)}}\right)f\left(
W_{\tau_{B(z_0,\frac{r}2)}}\right) \right],
$$
Thus by Lemma \ref{l:dis},
$$
f(x)\,\le\,  c\, M_1^2\,\E_x\left[f\left(
W_{\tau_{B(z_0,\frac{r}2)}}\right) \right] \,=:\, c\, M_1^2\, g(x)
$$
for some constant  $c$.
Since $g$ is harmonic for $W$ in $B(z_0,\frac{r}2)$, by the
Harnack inequality
for Brownian motion (for example, see \cite{B}),
$$
\frac1{c_1} g(y) \,\le\, g(x) \,\le\,c_1  g(y),   \quad x,y \in
B(z_0, \frac{r}4)
$$
for some constant $c_1 > 0$.
Thus by applying Theorem \ref{t:7.5} and Lemma \ref{l:dis} again,
we have that for every $ x,y \in B(z_0, \frac{r}4)$,
$$
f(x) \,\le\, c\, c_1\, M_1^2 \E_y\left[f\left(
W_{\tau_{B(z_0,\frac{r}2)}}\right) \right]
\,\le\, c^2 \,c_1\, M_1^4 \E_x\left[f\left(
\wh X^D_{\wh\tau_{B(z_0,\frac{r}2)}}\right);
\wh\tau_{B(z_0,\frac{r}2)}<\wh\zeta \right] \,=\,
c^2\, c_1\, M_1^4\, f(y).
$$
\qed

\medskip

It is well-known that the scale invariant Harnack inequality
implies the H\"{o}lder continuity of harmonic function
(for example, see section
2.3.2 of \cite{SC}).

\medskip

\begin{corollary}\label{c:Ho_C}
Every harmonic function for $\wh X^D$ is H\"{o}lder continuous.
\end{corollary}

\section{Martin representation in arbitrary bounded domains}

In this section we assume that
$D$ is an arbitrary bounded domain.
From Proposition 2.1, Theorem \ref{dualpr} and Proposition
\ref{s_feller_kd}, we see
that both $\wh X^D$ and $X^D$
satisfy the conditions (a)-(e) on page 560-561 in \cite{CK2}.
Thus by Theorem 3 in \cite{KW}
we can define
the Martin boundary for $\wh X^D$. In fact, we have a  stronger result.
We will state here for  $\wh X^D$ and $X^D$
 simultaneously. From now on, we fix a point $x_0\in D$ throughout
this paper.

Define
$$
M_D(x,y):=\left\{\begin{array}{ll}
\frac {\overline{G}_D(x,y)}{\overline{G}_D(x_0 , y)} \,= \,
\frac {G_D(x,y)}{G_D(x_0 , y)}     &\mbox{ if
$x \in D$ and $y \in D \setminus \{x_0\}$ }\\
1_{\{x_0\}}(x) & \mbox{ if  $y=x_0$  }
\end{array}\right.
$$
and
$$
\wh M_D(x,y):=\left\{\begin{array}{ll}
    \frac {\overline{G}_D(y,x)}{\overline{G}_D(y, x_0)} \,= \,
\frac {h_D(x)G_D(x,y)}{h_D(x_0)G_D(x_0 , y)}     &\mbox{ if
$x \in D$ and $y \in D \setminus \{x_0\}$ }\\
 1_{\{x_0\}}(x) & \mbox{ if  $y=x_0$  }
\end{array}\right.
$$

By (A4) and (A6), we know that for each $y \in D \setminus \{x_0 \}$
and $\eps>0$,
$M_D(\,\cdot\,,y)$ ($\wh M_D(\,\cdot\,,y)$ respectively)
is a harmonic function with respect to $X^D$ ($\wh X^D$ respectively)
in $D \setminus B(y, \eps)$
and for every $x \in D \setminus B(y, \eps)$
\begin{equation}\label{M_M}
M_D(x,y) = \E_x\left[M_D( X^D_{\tau_{D \setminus B(y, \eps)}}, y)
\right]
 \quad
\mbox{and}
\quad
\wh M_D(x,y) = \E_x\left[\wh M_D(\wh X^D_{\wh\tau_{D \setminus
B(y, \eps)}}, y); \wh\tau_{D \setminus
B(y, \eps)}<\wh\zeta\right].
\end{equation}

Using the Riesz decomposition theorem (Theorem \ref{t:Riesz}
and Theorem \ref{t:Riesz_dual}),
 Proposition \ref{p:ph} (with $U=D$),
the Harnack inequality (Corollary 5.8
in \cite{KS} and Theorem \ref{t:HP}) and
the H\"older continuity of harmonic functions (Theorem 5.5 in \cite{KS}
and Corollary \ref{c:Ho_C}),
one can follow the arguments in \cite{Ma} (see also Section 2.7 of
\cite{B} or \cite{Ta})
to show that the process $X^D$  ($\wh X^D$ respectively) has a Martin boundary
$\partial_MD$ ($\wh \partial_MD$ respectively) satisfying the
following properties.

\begin{description}
\item{(M1)} $D\cup \partial_M D$ and $D\cup \wh \partial_M D$ are
compact metric spaces;
\item{(M2)} $D$ is open and dense in $D\cup \partial_M D$ and in
$D\cup \wh \partial_M D$
and its relative topology coincides with its original topology;
\item{(M3)}  $M_D(x ,\, \cdot\,)$ can be extended  to $\partial_M$ uniquely
in such a way that,  $ M_D(x, y) $ converges to $M_D(x, w)$ as
$y\to w \in \partial_M D$,
the function $M_D(x, w)$  is jointly
continuous on $D\times \partial_M D$, and
$M_D(\cdot,w_1)\not= M_D(\cdot, w_2)$
if $w_1 \not= w_2$;
\item{(M4)}  $\wh M_D(x ,\, \cdot\,)$ can be extended  to
$\wh \partial_M$ uniquely
in such a way that,  $\wh M_D(x, y) $ converges to $\wh M_D(x, w)$ as
$y\to w \in \wh \partial_M D$,
the function $\wh M_D(x, w)$  is jointly
continuous on $D\times \wh \partial_M D$, and
$\wh M_D(\cdot,w_1)\not= \wh M_D(\cdot, w_2)$
if $w_1 \not= w_2$;
\end{description}

By repeating the argument in the proof of Proposition 5.1 in \cite{KS1},
we have the following.

\medskip

\begin{prop}\label{M_h_dual}
For every  $w\in   \partial_M D$ ($w\in \wh  \partial_M D$ respectively),
$x\mapsto M_D(x, w)$  ($x\mapsto
\wh M_D(x, w)$ respectively) is  harmonic with respect to $ X^D$
($\wh X^D$ respectively).
\end{prop}

\pf
We include the proof here for $\wh X^D$ for the reader's convenience.
Fix $ w\in  \partial_M D$ and  a relatively compact open sets
$U \subset \overline{U}
\subset U_1 \subset \overline {U_1}$ in $D$.
Let $\delta:= \frac12$dist$(U, \partial U_1)$.
Choose a sequence $\{y_n\}_{n \ge 1}$ in $D\setminus
\overline{U_1}$ converging
 to $w$ in  $D \cup \partial_M D$ so that
$$
\wh M_D(x, w)= \lim_{n \to \infty} \wh M_D(x, y_n).
$$
Since $\wh M_D(\,\cdot\,,y_n)$ is harmonic in a neighborhood of $U$
for every $n \ge 1$,
we have
$$
 \E_x\left[\wh M_D(\wh X^D_{\tau_{U}}, y_n)\right]\,=\,
\wh M_D(x, y_n), \quad x \in U.
$$
Using the
 Harnack inequality (Theorem \ref{t:HP}),
we have for every $z \in \partial U$,
$$
\wh M_D(z, y_n) =\frac{\overline{G}_D(y_n,z)}{\overline{G}_D(y_n, x_0)}
\,\le\, c_1 \frac{\overline{G}_D(y_n,x_0)}{\overline{G}_D(y_n, x_0)}
\,=\, c_1, \qquad n \ge 1,
$$
for some $c_1=c_1(\delta,D)>0$.
Thus
by the bounded convergence theorem,
$$
\lim_{n \to \infty} \E_x\left[\wh M_D(\wh X^D_{\wh\tau_{U}}; y_n);
\wh\tau_{U}<\wh\zeta\right]\,=\,
 \E_x\left[\wh  M_D(\wh  X^D_{\wh\tau_{U}}, w);
\wh\tau_{U}<\wh\zeta\right]\,=
\, \wh M_D(x, w), \quad x \in U.
$$
\qed

\medskip

Recall that a positive harmonic function $u$ with respect to $ X^D$
($\wh X^D$ respectively) is
said to be minimal if $v$ is positive harmonic with respect to $ X^D$
($\wh X^D$ respectively)
and $v\le u$ imply that $v$ is a constant multiple of $u$. The minimal
Martin boundaries of $ X^D$ and $\wh X^D$ are defined as
$$
\partial_mD=\{z\in \partial_MD: M_D(\cdot, z) \mbox{ is
minimal harmonic with respect to $X^D$ }\}
$$
and
$$
\wh \partial_mD=\{z\in \wh  \partial_MD:\wh  M_D(\cdot, z) \mbox{ is
minimal harmonic with respect to $\wh X^D$ }\}
$$
respectively.
Since $M_D(x_0, y) = 1$ for every $ y \in (D \cup \partial_M D)
\setminus \{x_0 \}$,
using the Harnack inequality and
the H\"older continuity of harmonic functions,
we can show that, for any compact subset $K$ of $D$, the family
$\{M_D(\cdot, w): w\in \partial _MD\}$ and $\{\wh M_D(\cdot, w):
w\in \partial _MD\}$
are uniformly bounded and
equicontinuous on $K$.
One can then apply the Ascoli-Arzel\`{a} theorem to prove
that, for every excessive function $f$ of $X^D$,
there exist a unique Radon measure $\nu_1$
on $D$ and a unique finite measure $\nu_2$ on $\partial_m D$
such that
\begin{equation}\label{eqn:rep}
 f(x)=\int_D G_D(x, y) \nu_1(dy)+
\int_{\partial_m D} M_D(x, z) \nu_2(dz),
\end{equation}
and $f$ is harmonic in $D$ with respect to $X$ if and only
if $\nu_1=0$.
Similarly,
for every excessive function $f$ of $\wh X^D$,
there exist a unique Radon measure $\mu_1$
on $D$ and a unique finite measure $\mu_2$ on $\wh \partial_m D$
such that
\begin{equation}\label{eqn:6.2_dual}
 f(x)=\int_D \overline{G}_D(y, x) \mu_1(dy)+
\int_{\wh \partial_m D} \wh M_D(x, z) \mu_2(dz),
\end{equation}
and $f$ is harmonic with respect to $\wh X^D$ if and only
if $\mu_1=0$. (See Section 2.7 of \cite{B})

\section{Martin boundary and Boundary Harnack principle for
$\wh X^D$}

In this section, we will, under some assumption on the domain,
identify the Martin boundary of the dual process with the
Euclidean boundary and prove a boundary Harnack principle for
the dual process.

Recall that a bounded domain $D$ is said to be Lipschitz
if there is a localization radius
$r_0>0$  and a constant
$\Lambda >0$
such that
for every $Q\in \partial D$, there is a
Lipschitz  function
$\phi_Q: \R^{d-1}\to \R$ satisfying
$| \phi_Q (x)- \phi_Q (z)| \leq \Lambda
|x-z|$, and an orthonormal coordinate
system $CS_Q$ with origin at $Q$ such that
$$
B(Q, r_0)\cap D=B(Q, r_0)\cap \{ y=(y_1, \cdots, y_{d-1}, y_d)=:
(\tilde y, y_d)\mbox{ in } CS_Q: y_d > \phi_Q (\tilde y) \}.
$$
The pair $(r_0, \Lambda)$ is called the
characteristics of the Lipschitz domain $D$.

We first recall the scale invariant boundary Harnack principle for $ X^D$ in
bounded Lipschitz domains from \cite{KS1}.

\medskip

\begin{thm}\label{BHP2}  (Theorem 4.6 in \cite{KS1})
Suppose $D$ is a bounded Lipschitz domain. Then there exist  constants
$M_2, c >1$  and $r_2 >0 $, depending on
$\mu$ only via the rate at which $\max_{1 \le i \le d}$
goes to zero such that for every $Q \in \partial D$, $r < r_2$ and any
nonnegative functions
$u$ and $v$ which are harmonic with respect to $X^D$ in $D \cap B(Q, M_2r)$
and vanish continuously on
$\partial D  \cap B(Q, M_2r)$, we have
\begin{equation}\label{e:BHP0}
\frac{ u(x)}{ v(x)}
\, \le \,c\, \frac{ u(y)}{ v(y)}  \quad  \mbox{ for any }
x,y \in D \cap B(Q, r).
\end{equation}
\end{thm}

\medskip

It is well-known  that for diffusions, the scale-invariant boundary
Harnack principle can be used to
prove  the H\"{o}lder continuity
of the ratio of two harmonic functions vanishing continuously near
the boundary. We omit the
proof of the next lemma. The proof can be found in \cite{B}
(see \cite{Bo} for the extension to jump
processes).

\medskip

\begin{lemma}\label{l:Ho}
Suppose $D$ is a bounded Lipschitz domain. Then there exist positive
constants $r_2$, $M_2$,
$C$ and $a$ depending on $D$ such
that for any $Q \in \partial D $, $r < r_2$ and nonnegative
functions $u, v$ which are harmonic with respect to $X^D$
in $D \cap B(Q, M_2 r)$ and vanish continuously on
$\partial D  \cap B(Q, M_2r)$, the limit $\lim_{D \ni x\to w}u(x)/v(x)$
exists for every $w\in \partial  D \cap B(Q, r)$.
\end{lemma}

\medskip

In this section we consider two bounded domains $U$ and $D$ with
$U \subset D$. We will not exclude the case $U=D$.
Let $x_U\in U$ (if $U=D$, $x_U=x_0$) and
define
$$
\wh M_{D,U}(x,y):=\left\{\begin{array}{ll}
\frac{h_D(x_U) G_U(y,x)}{h_D(x) G_U(y,x_U)}    &\mbox{ if
$x \in U$ and $y \in U \setminus \{x_U\}$ }\\
1_{\{x_U\}}(x) & \mbox{ if  $y=x_U$  }
\end{array}\right.
$$
Note that $\wh M_{D,U}(x,y)=\wh M_{D}(x,y)$ if $D=U$.
Similar to the proof of Proposition \ref{s_feller},
one can show that  $\wh X^{D,U}$ (with $X^U$ as a dual process)
satisfies the conditions (a)-(e) on page 560-561 in \cite{CK2}.
Thus by Theorem 3 in \cite{KW}, we
can define the Martin boundary $\wh \partial_M U$
for the process $\wh X^{D,U}$.
Moreover, one can prove that for every  $w\in \wh  \partial_M D$,
$x\mapsto
\wh M_{D,U}(x, w)$   is  harmonic with respect to $\wh X^{D}$ in $U$ using
an argument  similar to that of the proof of Proposition \ref{M_h_dual}.
Let $\wh \partial_m U$ be the minimal Martin boundary of $\wh X^{D,U}$.
 We also have the Martin representation:
for every harmonic function $f$ of $\wh X^D$ in $U$,
there is  a unique finite measure $\mu_1$ on $\wh \partial_m U$
such that
\begin{equation}\label{eqn:6.2_dual_kill}
f(x)= \int_{\wh \partial_m U} \wh M_{D,U}(x, z) \mu_1(dz).
\end{equation}

Suppose $U$ is a bounded Lipschitz domain.
We observe that for $y\not=x_U$,
\begin{equation}\label{e:M_M_kill}
\wh M_{D,U}(x,y) \,=\,
\frac {h_D(x_U)G_U(y,x)}{h_D(x)G_U(y, x_U)}.
\end{equation}
$G_U(\,\cdot\,,x)$ and $G_U(\,\cdot\,, x_U)$ are
harmonic with respect to $X^U$ near the boundary.
Moreover they  vanish continuously on  the boundary by Theorem 2.6
in \cite{KS1}. Thus
from Lemma \ref{l:Ho}, we immediately get the following

\medskip

\begin{prop}\label{p:M_conti_kill}
Suppose $U$ is a bounded Lipschitz domain. Then
$ \wh M_{D,U}(x, y) $ converges   as
$y\to w \in \partial U$.
\end{prop}

\medskip

The proposition above says that the Martin boundary is a subset of
$\partial U$.
We write the limit above as $ \wh M_{D,U}(x, w) $  for $(x,w) \in U
\times \partial U$.
Let
$N_{U}(x,w):=  \lim_{U \ni y \to w} \frac {G_U(y,x)}{G_U(y, x_U)}$
so that
$$
\wh M_{D,U}(x, w) =
\frac {h_D(x_U) N_{U}(x,w)   }{h_D(x)}.
$$

We will show that
the (minimal) Martin boundary $\wh \partial_m U$
with respect to $\wh X^{D,U}$
coincides with
the Euclidean boundary if $D$ and $U$ are bounded
$C^{1,1}$. Let $\rho_U(x)$ be the distance between $x$ and $\partial U$.

\medskip

\begin{thm}\label{t:dm_est_kill}
Suppose $D$ and $U$ are bounded $C^{1,1}$ domains with $U \subset D$.
Then there exists a constant
$c=c(D,U)$ such that
\begin{equation}\label{e:dm_est_kill}
\frac1{c}\, \frac{\rho_U(x)}{|x-w|^{d}} \,\le
\,\wh M_{D,U}(x, w)  \,\le\, c\,\frac{1}{|x-w|^{d}}
\end{equation}
and
\begin{equation}\label{e:dm_est}
\frac1{c}\, \frac{1}{|x-w|^{d}} \,\le
\,\wh M_D(x, w)  \,\le\, c\, \frac{1}{|x-w|^{d}}
\end{equation}
\end{thm}

\pf
By the Green function estimates for $X^U$ (Theorem 6.2 in \cite{KS}), we have
$$
c_1 \,\frac{\rho_U(x)}{|x-w|^{d}}\, \le\, N_{U}(x,w)\,
\le\, c_1\, \frac{\rho_U(x)}{|x-w|^{d}}
$$
for some positive constants $c_1$ and $c_2$.
Thus by (\ref{e:rho}) and the fact that $ \rho_U(x) \le \rho(x) \le $
diam$(D)$
$$
\frac1{c}
 \, \frac{h_D(x_U)\rho_U(x) }{|x-w|^{d}} \,\le
\,\wh M_{D,U}(x, w)  \,\le\, c\,
 \frac{h_D(x_U)}{|x-w|^{d}}.
$$
\qed

\medskip

The above  implies

\begin{prop}\label{l:dis2}If $D$ and $U$ are  bounded $C^{1,1}$
domains with $U \subset D$,
 $\wh M_{D,U}(\,\cdot\,, w_1)  \not= \wh M_{D,U}(\,\cdot\,, w_2)$
if $  w_1 \not= w_2$.
 \end{prop}

\medskip

Moreover, one can follow the argument in the proof of Theorem 4.4
\cite{CK} and show that
$\wh M_{D,U} (x, w)$ is minimal harmonic.
Thus the minimal Martin boundary
of $\wh X^{D,U}$ is the same as
the Euclidean boundary in the case when $D$ and $U$ are
bounded $C^{1,1}$ domains with $U \subset D$.

\medskip

\begin{thm}\label{t:MB3}
Assume that either $D$ and $U$ are bounded $C^{1,1}$ domains with
$U \subset D$
or $U$ is bounded Lipschitz domain
with $ \overline{U} \subset D$.   Then
for every harmonic function $f$ of $\wh X^D$ in $U$,
there is  a unique finite measure $\mu_1$ on $\partial U$
such that
\begin{equation}\label{rep_kill}
 f(x)= \int_{\partial U} \wh M_{D,U}(x, z) \mu_1(dz),
\end{equation}
$\partial D$.
\end{thm}

\pf The case when $D$ and $U$ are
bounded $C^{1,1}$ domains with $U \subset D$
has already been dealt with in the paragraph before the theorem.
In the case when  $U$ is bounded Lipschitz domain
with $ \overline{U} \subset D$, $\wh M_{D,U}(x, z)$ is comparable to
$N_U (x, z)$.
One can easily modify the argument in page 193-194 of \cite{B} to prove
the theorem.
We omit the details.
\qed

\medskip

Now we are ready to prove the boundary Harnack principle for $\wh X^D$.
If $D$ is a bounded $C^{1,1}$ domain, then
it is easy to check that there
exists $R>0$ such that for any $x\in \partial D$ and $r\in (0, R)$,
$B(x, r)\cap D$ is connected.

\medskip

\begin{thm}\label{BHP1} (Boundary Harnack principle)
Suppose $D$ be a bounded $C^{1,1}$ domain in $\R^d$ and $R$ is the constant
above.
Then for any $r\in (0, R)$ and $z_0\in \partial D$,
there exists a constant $c>1$
such that for any nonnegative harmonic functions $u, v$
in $D \cap B(z_0, r)$ with respect to $\wh X^D$ with $uh_D$ and
$vh_D$ vanishing
 continuously on
$\partial D \cap B(z_0, r)$, we have
$$
\frac{ u(x)}{ v(x)}
\, \le \,c\, \frac{ u(y)}{ v(y)}  \quad  \mbox{ for any }
x,y \in D\cap B(z_0, r/2).
$$
\end{thm}

\pf
One can find a bounded $C^{1,1}$ domain $U$ such that
$D \cap  B(z_0, 2r/3)\subset U \subset \overline{U}
\subset D \cap B(z_0, r).$
Fix $x_U \in U$  and  let
$$
M_2(x, z)
:= \wh M_{D,U}(x, z)= \lim_{U \ni y  \to z} \frac{h_D(x_U) G_{U}(y,x)}
{h_D(x) G_{U}(y,x_U)}.
$$
Since $u, v$ are harmonic in $U$ with respect to $\wh X^D$,
by Theorem \ref{t:MB3}, there exist finite measures $\mu_1$
and $\nu_1$ on $\partial U$
such that
$$
u(x)= \int_{\partial U} M_{2}(x, z) \mu_1(dz) \quad \mbox{and}\quad
v(x)=\int_{\partial U} M_{2}(x, z) \nu_1(dz), \quad x \in U.
$$
Let
$N_2(x,z):=  \lim_{U \ni y \to z} \frac {G_{U}(y,x)}{G_{U}(y, x_U)}$
so that
$$
M_2(x, z) =
\frac {h_D(x_U) N_2(x,z)   }{h_D(x)}.
$$
Let $G^{0}_{U}$ be the Green function of the
Brownian motion $W$ in $U$.
Define the Martin kernel $M_1(x, z)$ for the Brownian motion $W$ in $U$:
$$
M_1(x,z):=  \lim_{U \ni y \to z} \frac {G^0_{U}(x,y)}{G^0_{U}(x_U,y)}.
$$
Since $U$ is bounded $C^{1,1}$, by Theorem 7.7 in \cite{KS},
there exists a constant
$c_1=c_1(x_U, U)$ such that
\begin{equation}\label{e:comp}
\frac1{c_1}\, M_1(x,z) \,\le\,N_2(x,z) \,\le\,c_1\, M_1(x,z).
\end{equation}
Let
$$
u_1(x):= \int_{\partial U}  M_{1}(x, z) \mu_1(dz) \quad \mbox{and}\quad
v_1(x):=\int_{\partial U}  M_{1}(x, z) \nu_1(dz), \quad x \in U.
$$
By (\ref{e:comp}), we have for every $x \in U$
\begin{eqnarray*}
&&\frac{u(x)}{v(x)}\,=\,\frac{\int_{\partial U} M_{2}(x, z) \mu_1(dz)}
{\int_{\partial U} M_{2}(x, z) \nu_1(dz)}
 \,=\,\frac{\int_{\partial U} N_{2}(x, z) \mu_1(dz)}
{\int_{\partial U} N_{2}(x, z) \nu_1(dz)}\\
 &&\le\,c_1^2\frac{\int_{\partial U} M_{1}(x, z) \mu_1(dz)}
{\int_{\partial U} M_{1}(x, z) \nu_1(dz)}
\,=\,c_1^2\frac{u_1(x)}{v_1(x)}\,\le\, c_1^4 \frac{u(x)}{v(x)}.
\end{eqnarray*}

Since $u_1, v_1$ are harmonic for the Brownian motion $W$ in $U$ and vanish
continuously on $\partial U \cap \partial D$, by the boundary Harnack
principle
for Brownian motion (for example, see \cite{B}),
$$
\frac{u_1(x)}{v_1(x)} \,\le\, c_2 \frac{u_1(y)}{v_1(y)}, \quad x,y \in
D \cap B(z_0, \frac{r}2)
$$
for some constant $c_2 > 0$.
Thus for every $ x,y \in
D \cap B(z_0, \frac{r}2)$
$$
\frac{u(x)}{v(x)} \,\le\,c_1^2\frac{u_1(x)}{v_1(x)}\,\le\,c_2 c_1^2
 \frac{u_1(y)}{v_1(y)}\,\le\,c_2 c_1^4 \frac{u(y)}{v(y)}.
$$
\qed

\section{Schr\"{o}dinger operator
in arbitrary bounded domains}

In this section we discuss the Schr\"{o}dinger operator
in arbitrary bounded domains.
Using results in the previous sections, one can check that
 $X^D$ and $\wh X^D$ satisfy
the condition (a)-(f) and (6.1) in \cite{CK2} with the reference
measure $\xi_D$. Thus the main results in \cite{CK2} are true for $X^D$
and $\wh X^D$ with the reference measure $\xi_D$. In this section
we will use the main results in \cite{CK2} and
 state carefully for $X^D$  with respect to the Lebesgue measure.

Recall that a measure $\nu$ on $D$ is said to be a smooth measure
of $X^D$ with respect to the reference measure $\xi_D$
if there is a positive continuous additive functional
(PCAF in abbreviation) $A$ of $X^D$ such that for all
bounded nonnegative function $f$ on $D$,
\begin{equation}\label{eqn:Revuz1}
\int_D f(x) \nu (dx) = \lim_{t\downarrow 0}
\E_{\xi_D} \left[ \frac1t \int_0^t f(X^D_s) dA_s \right].
\end{equation}
The additive functional $A$ is called the positive continuous
additive functional of $X^D$ with Revuz measure $\nu$
with the reference measure $\xi_D$.
It is known (see \cite{Re}) that for any $x\in D$, $\alpha\ge 0$ and
bounded nonnegative function $f$ on $D$,
$$
\E_x \int^{\infty}_0e^{-\alpha t}f(X^D_t) dA_t =
\int^{\infty}_0e^{-\alpha t}\int_D
\overline{q}^D(t, x, y)f(y) \nu (dy)dt,
$$
and thus we have any $x\in D$,  $t>0$ and
bounded nonnegative function $f$ on $D$,
$$
\E_x \int^t_0f(X^D_s) dA_s=\int^t_0\int_D
\overline{q}^D(s, x, y)f(y) \nu (dy)ds.
$$
Therefore by the monotone convergence theorem we
have any $x\in D$,  $t>0$ and
nonnegative function $f$ on $D$,
\begin{equation}\label{eqn:revuz2}
\E_x \int^t_0f(X^D_s) dA_s=\int^t_0\int_D
\overline{q}^D(s, x, y)f(y) \nu (dy)ds.
\end{equation}
For a signed measure $\nu$, we use $\nu^+$ and $\nu^-$
to denote its positive and negative parts respectively.
If $\nu^+$ and $\nu^{-}$ are smooth measures of $X^D$ with respect to
the reference measure $\xi_D$ and
$A^{+}$ and $A^{-}$ are PCAFs of $X^D$ with Revuz measures
$\nu^+$ and $\nu^{-}$ respectively with respect to the reference
measure $\xi_D$, then we call
$A:=A^+-A^{-}$ of $X^D$ the continuous additive functional
of $X^D$ with (signed) Revuz measure $\nu$ with respect
to the reference measure $\xi_D$.

A measure $\eta$ on $D$ is said to be a smooth measure
of $X^D$ with respect to the Lebesgue measure if $h_D\eta$
is a smooth measure of $X^D$ with respect to the reference
measure $\xi_D$. From now on, whenever we speak of a smooth
measure of $X^D$, we mean a smooth measure of $X^D$ with
respect to the Lebesgue measure unless explicitly mentioned otherwise.
If $\eta$ is a smooth measure of $X^D$ and $A$ is
the PCAF of $X^D$ with Revuz measure $h_D\eta$
with respect to the reference measure $\xi_D$, then by
(\ref{eqn:revuz2}), we have any $x\in D$,  $t>0$ and
bounded nonnegative function $f$ on $D$,
\begin{equation}\label{eqn:revuz3}
\E_x \int^t_0f(X^D_s) dA_s=\int^t_0\int_D
q^D(s, x, y)f(y) \eta (dy)ds.
\end{equation}
The additive functional $A$ in the equation above
is called the PCAF of $X^D$ with Revuz measure $\eta$
with respect to the Lebesgue measure. From now on,
whenever we speak of an additive functional with a
given Revuz measure, we mean an additive functional with a
given Revuz measure with respect to the Lebesgue
measure unless explicitly mentioned otherwise.

Recall $\K_{d, 2}$ from the Definition \ref{d:kc}.

\medskip

\begin{prop}\label{Katoissm}
Any measure $\nu$ in $\K_{d, 2}$ is a smooth measure
of $X^D$ with respect to the Lebesgue measure, or equivalently,
$h_D\nu$ is a smooth measure
of $X^D$ with respect to $\xi_D$.
\end{prop}

\pf By the definition of $\K_{d, 2}$ we can easily check
that the function
$$
\overline{G}_D(h_D\nu)(x)=G_D\nu(x)
$$
is bounded continuous in $D$. Thus, by Definition
IV.3.2 of \cite{BG}, $\overline{G}_D(h_D\nu)$ is a
regular potential. Moreover
$X^{D}$ and $\wh X^{D}$ are Hunt processes with the
strong Feller property
and they are in
the strong duality with respect to $\xi_D$
(Propositions \ref{dualpr} and \ref{s_feller_kd}).
Consequently $h_D\nu$ charges
no semi-polar set by Theorem VI.3.5 of \cite{BG}.
Now we can apply Theorem VI.1 of \cite{Re} to
conclude that $h_D\nu$ is a smooth measure
of $X^D$ with respect to $\xi_D$, or equivalently,
$\nu$ is a smooth measure
of $X^D$ with respect to the Lebesgue measure.
\qed

\medskip

If $\nu$ is a signed measure on $D$ such that $\nu^+$
and $\nu^-$ are smooth measures of $X^D$ and if
$A^+$ and $A^-$ are the PCAF of $X^D$ with Revuz measures
$\nu^+$ and $\nu^-$ respectively, we call
$A:=A^+-A^{-}$ of $X^D$ the continuous additive functional
of $X^D$ with (signed) Revuz measure $\nu$.

We recall the definitions of the
class of measures from
\cite{C2} and \cite{CS6} and specify it
for $X^D$ with the reference measure $\xi_D$ in an
bounded domain $D$. We also give a definition for a class of smooth
measures with respect to the Lebesgue measure.
In the following,
$d$ denotes the diagonal of $D\times D$.

\medskip

\begin{defn}\label{df:3.1}
Let $\nu$ be a signed smooth measure of $X^D$ with
respect to $\xi_D$ and
define $|\nu|:=\nu^++\nu^-$. $\nu$
is said to be in the class $\S^{\xi_D}_\infty (X^D)$
if for any $\eps>0$ there is a Borel
subset $K=K(\eps)$ of finite $| \nu|$-measure
and a constant $\delta = \delta (\eps) >0 $ such that
for all $(x, z)\in (D\times D)\setminus d$,
\begin{equation}\label{eqn:3.A}
\int_{D\setminus K}
\frac{\overline{G}_D(x, y)\overline{G}_D(y, z) }{\overline{G}_D(x, z)}
\, |\nu | (dy)=\int_{D\setminus K}
\frac{G_D(x, y){G}_D(y, z) }{{G}_D(x, z)} \,
\frac{|\nu | (dy)}{h_D(y)}\le\eps
\end{equation}
and for all measurable set $B\subset K$ with $| \nu |(B)<\delta$
all $(x, z)\in (D\times D)\setminus d$,
\begin{equation}\label{eqn:3.B}
 \int_B \frac{\overline{G}_D(x, y)\overline{G}_D(y, z) }
{\overline{G}_D(x, z)} \, |\nu |(dy)
=\int_B \frac{{G}_D(x, y){G}_D(y, z) }{{G}_D(x, z)} \,
\frac{|\nu |(dy)}{h_D(y)} \le\eps.
\end{equation}
A function $q$ is said to be in the class
$\S^{\xi_D}_\infty (X^D)$
if $\nu(dx) :=q(x)\xi_D(dx)$ is in the corresponding space.

A signed smooth $\nu$ of $X^D$
is said to be in the class $\S_{\infty}(X^D)$ if
$h_D(x)\nu(dx)$ is in the class $\S^{\xi_D}_\infty (X^D)$.
\end{defn}

\medskip

For a continuous additive functional $A$ of $X^D$ with Revuz measure
$\nu$, we define $e_{A}(t)=\exp( A_t),$ for $t\geq 0.$
For $y\in D$, let $X^{D,y}$ denote the $h$-conditioned process obtained
from  $X^D$ with
$h(\cdot)= G_D(\cdot , y)$ and let $\E_x^y$ denote the expectation
for $X^{D,y}$ starting from $x\in D$. We will use $\tau^y_D$ to denote
the lifetime of the process $X^{D,y}$.

In the remainder of this section, we assume that
$\nu\in \S_\infty(X^D)$
and $A$ is the CAF of $X^D$ with Revuz measure $\nu$
(with respect to the Lebesgue measure).
Note that  $A$ is also the CAF of $X^D$ with Revuz measure $h_D(x)\nu(dx)$
with respect to the reference measure $\xi_D$.

The CAF $A$ gives rise to a Schr\"odinger semigroup
$$ Q^D_t f(x):= \E_x \left[ e_{A}(t)f(X^D_t) \right].
$$

The function $x\mapsto \E_x [e_{A}(\tau_D)]$ is called the gauge function
of $\nu$. We say $\nu$ is {\it gaugeable}
if $ \E_x [ e_{A} (\tau_D) ]$ is finite
for some $x\in D$. From now on we will assume that
$\nu$ is gaugeable. It follows from \cite{C2} and \cite{CS6}
that the gauge function $x\mapsto \E_x \left[
e_{A}(\tau_D)\right]$ is bounded on $D$.
Note that since $h_D(x)\nu(dx) \in \S^{\xi_D}_\infty (X^D)$,
it follows again from \cite{C2} and \cite{CS6} that
$$
\sup_{(x, y)\in (D\times D)\setminus d} \E_x^y \left[
|A|_{\tau_D^y}\right]< \infty
$$
(see \cite{C2}) and therefore by Jensen's inequality
\begin{equation}\label{eqn:jen}
\inf_{(x, y)\in (D\times D)\setminus d} \E_x^y [ e_{A}
(\tau_D^y)]>0.
\end{equation}

By Lemma 3.5 of \cite{C2}, the Green function for the
Schr\"odinger semigroup $\{Q^D_t, \, t\geq 0\}$ with respect to $\xi_D$ is
\begin{equation}\label{eqn:3.6}
\overline{V}_D(x, y)=\E_x^y
\left[ e_{A}(\tau_D^y) \right] \overline{G}_D(x, y),
\end{equation}
that is,
\begin{equation}\label{eqn:new3.8}
\int_D \overline{V}_D(x, y) f(y) \, \xi_D(dy) = \int_0^\infty
Q^D_t f (x) \, dt
= \E_x \left[ \int_0^\infty e_{A}(t) f(X^D_t) \, dt \right]
\end{equation}
for any Borel function $f\geq 0$ on $D$.
Let
$$
V_D(x, y) := \overline{V}_D(x, y) h_D(y)
$$
so that
\begin{equation}\label{eqn:new3.8_1}
\int_D V_D(x, y) f(y) \, dy = \int_0^\infty Q^D_t f (x) \, dt
= \E_x \left[ \int_0^\infty e_{A}(t) f(X^D_t) \, dt \right]
\end{equation}
Thus $V_D(x, y)$ is the Green function for the
Schr\"odinger semigroup $\{Q^D_t, \, t\geq 0\}$ with respect to the Lebesgue
measure. Since $G_D(x, y) = \overline{G}_D(x, y) h_D(y)$,
by Theorem 3.6 in \cite{C2} and (\ref{eqn:jen}),
$V_D(x, y)$ is comparable to $G_D(x, y)$ on $D\times D \setminus d$,
where $d$ denotes the diagonal of $D\times D$.
From Theorem 3.4 \cite{CK2} and the continuity of $h_D$, we see
that $V_D(x, y)$ is continuous on
$(D\times D)\setminus d$.

Let $u(x, y):=\E_x^y \left[ e_{A}(\tau_D^y) \right]$
for $y\in D$, and define $u(x, w):=
\E_x^w \left[ e_{A}(\tau_D^w) \right]$ for $w\in \partial D$,
where $\E_x^w$ is the expectation for
the conditional process of $X^D$ obtained through $h$-transform
with $h(\cdot )=M(\cdot, w)$. Recall that $\partial_m D$ is the
minimal Martin boundary
of $X^D$.

\medskip

\begin{thm}\label{T:C1}
The following properties hold.
\begin{description}
 \item{(1)} For $w\in \partial_m D$ and $x\in D$,
 $\lim_{D \ni y\to w} u(x, y)=u(x, w)$.
The conditional gauge function $u(x, w)$
 is jointly continuous on $D\times \partial_m D$.

\item{(2)} For every $x\in D$ and $w\in \partial_m D$, $K_D(x, w):=
  \lim_{D \ni y\to w} \frac{V_D(x, y)}{V_D(x_0, y)}$ exists and is finite.
Furthermore,
\begin{equation}\label{eqn:K1}
K_D(x, w)=M_D(x, w)\frac{u(x, w)}{u(x_0, w)}
\end{equation}
and so $K(x, w)$ is jointly continuous on $D\times \partial_m D$;

\item{(3)} Assume $D$ is a bounded $C^{1,1}$ domain and let
$\rho_D(x)$ be the distance between
$x$ and $ \partial D$, then there is a constant $c>1$ such that
\begin{equation}\label{eqn:K2}
 c^{-1}\frac{\rho_D(x)}{|x-w|^{d}}
\leq K_D(x, w) \leq c\, \frac{\rho_D(x)}{|x-w|^{d}}.
\end{equation}
\end{description}
\end{thm}

\pf (1) and (2) are proved in Theorem 3.4 in \cite{CK2}
(also see section 6 in \cite{CK2} for the extension).
Estimate (\ref{eqn:K2}) follows directly from (\ref{eqn:K1}) above
and Theorem 7.7 in \cite{KS}.
\qed

\medskip

\begin{defn}\label{D:3.5}
A  Borel function
$u$ defined on $U$ is said to be $\nu$-harmonic for $X^D$
in an open subset $U$ of $D$ if
$$
\E_x \left[ e_{A}(\tau_B) |u(X^D_{\tau_B})| \right]
<\infty \quad \mbox{and} \quad
\E_x \left[ e_{A}(\tau_B) u(X^D_{\tau_B})
\right] = u(x), \quad  x\in B,
$$
for every open set $B$ whose closure is a compact subset of $U$.
If $U=D$, then $u$ is said to be $\nu$-harmonic for $X^D$ .
\end{defn}

\medskip

The following two theorems are proved in \cite{CK2}. (Lemma 3.6,  Theorems
5.11-5.12, Theorem 5.14-5.16 \cite{CK2}. Also see section 6 in
\cite{CK2} for the extension.)

\medskip

\begin{thm}\label{T:5.15}
For
 every $z \in \partial_m D$, $x \mapsto K_D(x,z)$ is
a minimal  $\nu$-harmonic function
of $X^D$. That is, if $h$ is a $\nu$-harmonic
function of $X^D$ and $0\leq h(x) \leq K_D(x,z)$,
then $h(x)=c K_D(x,z)$ for some constant $c\leq 1$. Moreover,
for $z_1 \not= z_2 \in
\partial_m D$, $ K_D(\cdot,z_1) \not\equiv K_D(\cdot,z_2)$.
\end{thm}

\medskip

Recall that a nonnegative Borel  function
$f$ defined on $D$
is said to be
$\nu$-excessive for $X^D$ if for every $x\in D$ and $t>0$,
$Q^D_t f(x) \le f(x)$ and
$$
\lim_{t \downarrow 0} Q^D_t f(x) = f(x).
$$

\medskip

\begin{thm}\label{T:5.16}
The minimal Martin
boundary for the Schr\"odinger semigroup $\{Q^D_t, t\geq 0\}$
can be identified with the minimal Martin boundary
$\partial_m D$ of $X^D$.
Furthermore for every $\nu$-excessive function $f$ of $X^D$
that is not identically infinite,
there is a unique Radon measure $\mu_1$ on $D$
and a unique finite measure $\mu_2$ on $\partial D$ such that
\begin{equation}\label{eqn:5.8}
 f(x) = \int_D V_D(x, y)  \mu_1 (dy) + \int_{\partial_m D} K_D(x, z)
\mu_2 (dz).
\end{equation}
Function $f$ is $\nu$-harmonic for $X^D$ if and only $\mu_1=0$.

Conversely, if $\mu_1$ is a Radon measure in $D$ such that
$\int_D V_D(x, y)  \mu_1 (dy)$
is not identically infinite and $\mu_2$ is  finite
measure on $\partial D$, then the function $f$ given by
(\ref{eqn:5.8}) is a nonnegative
$\nu$-excessive function of $X^D$
that is not identically infinite.
\end{thm}

\medskip

Therefore we conclude that the minimal Martin boundary is stable
under Feynman-Kac perturbation
if $\nu\in \S_\infty(X^D)$
such  the gauge function $x\mapsto \E_x \left[
e_{A}(\tau_D)\right]$ is bounded.
Furthermore we see that there is a one-to-one correspondence between
the space of excessive functions
(the space of positive harmonic functions) of $X^D$
that are not identically infinite
and the space of $\nu$-excessive functions
(the space of positive $\nu$-harmonic functions, respectively)
of $X^D$ that are not identically infinite
through measures $\mu_1$ and $\mu_2$.

Since the Martin measure $\mu_2$ is finite and $K_D(x, z)$
is jointly continuous on $D \times \partial_m D$, we
have the continuity for $\nu$-harmonic functions of $X^D$.

\medskip

\begin{thm}\label{t:nu_cont} If
$u\geq 0$ is $\nu$-harmonic for $X^D$,
then $u$ is continuous in $D$.
\end{thm}

\medskip

In \cite{KS1}, we have shown that
for every  bounded Lipschitz domain $D$,
there is a one-to-one correspondence between
the minimal Martin boundary $\partial_m D$ for $X^D$
and the Euclidean boundary
$\partial D$ (see Theorem 5.7 in \cite{KS1}). Thus from
Theorem \ref{T:5.16}, we have

\medskip

\begin{thm}\label{T:S_Lip}
For every  bounded Lipschitz domain $D$,
 the (minimal) Martin
boundary for the Schr\"odinger semigroup $\{Q^D_t, t\geq 0\}$
 can be identified with the Euclidean boundary.
Furthermore, for every $\nu$-excessive function $f$ of $X^D$
that is not identically infinite, $\partial D$
there is a unique Radon measure $\mu_1$ on $D$
and a unique finite measure $\mu_2$ on $\partial D$ such that
\begin{equation}\label{eqn:5.9}
 f(x) = \int_D V_D(x, y)  \mu_1 (dy) + \int_{\partial D} K_D(x, z) \mu_2 (dz).
\end{equation}
Function $f$ is $\nu$-harmonic for $X^D$ if and only $\mu_1=0$.
\end{thm}

\medskip

\begin{remark}
{\rm In \cite{KS3}, by using the Green function estimates
and our Proposition \ref{Katoissm}, we show that,
in fact, $\K_{d, 2}$ is contained in $\S_\infty(X^D)$ if
$D$ is bounded Lipschitz.}
\end{remark}

\vspace{.5in}
\begin{singlespace}

\end{singlespace}
\end{doublespace}
\end{document}